\documentclass[a4paper,11pt,reqno,twoside]{amsart}
\usepackage{amsmath}
\usepackage{amsfonts}
\usepackage{amssymb}
\usepackage{amsthm}
\usepackage{verbatim}
\theoremstyle{plain}
\newtheorem{thm}{Theorem}[section]
\newtheorem{cor}[thm]{Corollary}
\newtheorem{lemma}[thm]{Lemma}
\newtheorem{prop}[thm]{Proposition}
\newtheorem{defn}[thm]{Definition}
\theoremstyle{definition}
\newtheorem{remark}[thm]{Remark}
\newtheorem{example}[thm]{Example}
\numberwithin{equation}{section}

\newcommand{\N}{\mathbb{N}}
\newcommand{\C}{\mathbb{C}}
\newcommand{\Z}{\mathbb{Z}}

\newcommand{\h}{\mathfrak{h}}
\newcommand{\hs}{\mathfrak{h}^\ast}
\newcommand{\rl}{\mathfrak{r}}
\newcommand{\rls}{\mathfrak{r}^\ast}
\newcommand{\W}{\mathcal{W}}
\newcommand{\A}{{\mathcal{A}}}
\newcommand{\B}{{\mathcal{B}}}
\newcommand{\U}{{\mathcal{U}}}
\newcommand{\E}{{\mathcal{E}}}

\newcommand{\id}{\mathrm{Id}}

\newcommand{\lr}{{\mathrm{lr}}}

\renewcommand{\det}{\mathrm{det}}
\newcommand{\Dh}{{D_{\hs}}}
\newcommand{\Dr}{{D_{\rls}}}
\newcommand{\Ih}{{I_{\hs}}}
\newcommand{\Mh}{{M_{\hs}}}
\newcommand{\Mr}{{M_{\rls}}}
\newcommand{\1}{\mathbf{1}}
\newcommand{\halg}{$\h$-algebra}
\newcommand{\hprealg}{$\h$-prealgebra}
\newcommand{\hcoalg}{$\h$-coalgebroid}
\newcommand{\hbialg}{$\h$-bialgebroid}
\newcommand{\hHopfalg}{$\h$-Hopf algebroid}
\newcommand{\hHopfster}{$\h$-Hopf $\ast$-algebroid}
\newcommand{\hspa}{$\h$-space}
\newcommand{\wtt}{\widetilde\otimes}

\newcommand{\opp}{{\mathrm{opp}}}
\newcommand{\cop}{{\mathrm{cop}}}
\newcommand{\da}{{\dagger}}
\newcommand{\al}{\alpha}
\newcommand{\be}{\beta}
\newcommand{\ga}{\gamma}

\newcommand{\de}{\delta}
\newcommand{\si}{\sigma}
\newcommand{\la}{\lambda}
\newcommand{\om}{\omega}
\newcommand{\De}{\Delta}
\newcommand{\ep}{\varepsilon}
\newcommand{\io}{\iota}

\newcounter{mylist}

\newenvironment{pair}[2]{ \langle#1, #2\rangle}

\begin{document}
\date{December 10, 2004}
\title
{Pairings and actions for dynamical quantum groups}
\author{Erik Koelink and Yvette van Norden}
\thanks{The second author is supported by Netherlands Organisation for Scientific
Research (NWO) under project number 613.006.572.}
\address{Technische Universiteit Delft, EWI-DIAM,
Postbus 5031, 2600 GA Delft, the Netherlands}
\email{h.t.koelink@ewi.tudelft.nl, y.vannorden@ewi.tudelft.nl}

\begin{abstract}Dynamical quantum groups
constructed from a FRST-construction using a
solution of the quantum dynamical Yang-Baxter
equation are equipped with a natural pairing.
The interplay of the pairing with $\ast$-structures,
corepresentations and dynamical representations
is studied, and natural left and right actions are introduced.
Explicit details for the elliptic $U(2)$ dynamical
quantum group are given, and the pairing is
calculated explicitly in terms of elliptic hypergeometric
functions. Dynamical analogues of spherical and
singular vectors for corepresentations are
introduced.
\end{abstract}

\maketitle

\section{Introduction}

Dynamical quantum groups have been introduced by
Etingof and Varchenko \cite{EtinV}, \cite{EtinV2}
in order to provide an algebraic framework,
in terms of \hHopfalg s (see \S \ref{ssec:algebraicnotions}),
for the
study of the quantum dynamical Yang-Baxter equation
\eqref{eq:QDYBE} similar to the relation of
quantum groups and (constant) solutions of the
quantum Yang-Baxter equation. The construction is
an analogue of the Faddeev-Reshetikhin-Sklyanin-Takhtajan
(FRST) construction, which we recall in \S
\ref{ssec:FRST-construction}.
The quantum dynamical Yang-Baxter
equation arises in a natural way from the construction
of correlation functions and corresponding fusion and
exchange matrices in the study of the
Knizhnik-Zamolodchikov-equations and its difference
analogue, see  \cite{EtinFK}, \cite{EtinS} and
references given there.
In particular, this construction
attaches an algebraic framework to the elliptic
$R$-matrix \eqref{eq:ellR} related to $\mathfrak{sl}(2)$
involving both a spectral and a dynamical
parameter.  The corresponding algebraic framework
has been studied in detail by Felder and Varchenko
\cite{FeldV}, which predates \cite{EtinV},
\cite{EtinV2}, and we call this the elliptic $U(2)$
dynamical quantum group and this is the main example
in this paper, see \S \ref{ssec:FRST-construction}.

Another well-studied example is the dynamical
quantum group associated to the rational $R$-matrix \eqref{eq:ratR}
for $\mathfrak{sl}(2)$ that can be
obtained by a limit transition from the elliptic
$R$-matrix. The corresponding algebraic structure
is essentially simpler, since the dependence on the
spectral parameter is removed.
In \cite{KoelR} its corepresentation theory
has been studied, and it turns out that there is
a direct link to special functions, in particular to
Askey-Wilson and $q$-Racah polynomials which
occur in the description of matrix elements of
irreducible corepresentations. These
results have been partly extended to the
elliptic $U(2)$ dynamical quantum group in \cite{KoelvNR},
linking corepresentations to so-called elliptic
hypergeometric series, originally introduced
by Frenkel and Turaev \cite{FrenT} and
which are recalled in \S \ref{ssec:ellpairingmatrixelt}.

An important ingredient for this paper is the fact that
dynamical quantum groups arising from a FRST-construction
for which the $R$-matrix is a solution to the quantum
dynamical Yang-Baxter equation is self-dual in the sense
that there is a natural pairing between the dynamical
quantum group and its coopposite dynamical quantum group, which is
compatible with the algebraic structures, but which in
general is not non-degenerate. The pairing is defined
in terms of the $R$-matrix by Rosengren \cite{Rose},
and the example of the dynamical
quantum group associated to the rational $R$-matrix
for $\mathfrak{sl}(2)$ is completely worked out in \cite{Rose}.
A similar duality has been given in the context of
weak Hopf algebras by Etingof and Nikshych \cite[\S 5]{EtinN}
related to the case of a dynamical quantum group at a root of $1$,
in which case the pairing is non-degenerate. We recall
Rosengren's definitions in \S \ref{sec:pairing}, and
we extend it to include the case of a spectral parameter
in the $R$-matrix and to include a $\ast$-structure.
In particular, we work out the details of the pairing for the
elliptic $U(2)$ dynamical quantum group.
The algebraic structures for dynamical quantum groups are
\hHopfalg s, which are recalled in
\S \ref{ssec:algebraicnotions}, and in general
the pairing between two \hHopfalg s induces a natural
action of one \hHopfalg\ on the other, extending the
results for ordinary Hopf algebras and which can
be thought of a generalisation of the action
$(X\cdot f)(g) =\frac{d}{dt}\big\vert_{t=0} f(g\exp(tX))$
for $X\in\mathfrak{g}$, $f$ a polynomial on the
Lie group $G$ with Lie algebra $\mathfrak{g}$.
This also gives the
opportunity to go from corepresentations of one
\hHopfalg\  to dynamical representations of the other
\hHopfalg . In the case of the elliptic $U(2)$ dynamical
quantum group we work out the details, and in particular
we calculate the pairing between matrix elements
of irreducible corepresentations in terms of elliptic
hypergeometric series. This then allows us to give
a dynamical quantum group derivation of the quantum
dynamical Yang-Baxter equation for the elliptic hypergeometric
series,
biorthogonality relations, the Bailey transform
and the Jackson sum for elliptic hypergeometric series
already obtained by Frenkel and Turaev \cite{FrenT},
see \cite[Ch. 11]{GaspR} for more information and
references. For this we rely on a previous paper
\cite{KoelvNR} in which a slightly different dynamical
quantum group theoretic derivation of
biorthogonality relations, the Bailey transform
and the Jackson sum for elliptic hypergeometric series
are given.

Starting with the work of Koornwinder \cite{Koor} we
also have the interpretation of Askey-Wilson polynomials
as spherical functions for irreducible corepresentations of
the standard quantum $SU(2)$ group using so-called
twisted primitive elements, see also \cite{KoelAAM}, \cite{NoumM}
for extensions to arbitrary matrix elements. This idea
with the appropriate replacement of twisted primitive elements
by co-ideals has turned out to be very fruitful, especially
for quantum analogues of compact symmetric spaces and the
corresponding spherical functions in terms of
Askey-Wilson and Macdonald-Koornwinder polynomials, see
Dijkhuizen and Noumi
\cite{DijkN}, Noumi \cite{Noum},
Noumi, Dijkhuizen and Sugitani \cite{NoumDS}, and Letzter \cite{Letz},
and for a non-compact quantum group
and the Askey-Wilson functions see \cite{KoelS}.

Since the interpretation of Askey-Wilson polynomials on the
standard quantum $SU(2)$ group and dynamical
quantum group associated to the rational $R$-matrix
for $\mathfrak{sl}(2)$ is similar, a natural link between these
algebras is to be expected and the precise relation has been given
by Stokman \cite{Stok}. A key ingredient is the twisted
coboundary element of Babelon, Bernard and Billey \cite{BabeBB}
in the form discovered by Rosengren \cite{RoseCM} as a
universal element
intertwining the standard Cartan element with the
twisted primitive element of Koornwinder. This element is
only known explicitly for the case $\mathfrak{sl}(2)$, see
Buffenoir and Roche \cite{BuffR} for conjectural forms
of the twisted coboundary element for other simple
$\mathfrak{g}$, in particular $\mathfrak{g}=\mathfrak{sl}(n)$.
One may conjecture that the harmonic analysis on
dynamical quantum analogues (for the rational $R$-matrices)
of compact symmetric spaces may shed light on this matter.
This is the motivation for the discussion in \S \ref{sec:singular}
in which we consider an application of the action to
the notion of singular and spherical vectors in corepresentations
for dynamical quantum groups that are equipped with a pairing.
The general definitions have to be worked out for explicit
examples, and we intend to do so in a future paper.

For compact quantum groups there is a natural analogue of the
Haar functional, and there is an analogue of Schur's orthogonality
relations. For the dynamical
quantum group associated to the rational $R$-matrix
for $\mathfrak{sl}(2)$ the Haar functional
was introduced using the analogue
of the Peter-Weyl theorem and the Clebsch-Gordan decomposition,
see \cite[\S 7]{KoelR}.
This method is not generally applicable to dynamical quantum
groups, and we expect that the pairing and the actions
defined in this paper can give rise to an alternative definition.
It remains open if the elliptic beta integral, see
e.g. \cite[Exerc. 11.29]{GaspR}, can be given
a dynamical quantum group theoretic interpretation in this way.
See B\"ohm and Szlach\'anyi \cite{BohS} for integrals in the
context of Takeuchi's $\times_R$-bialgebras with an antipode.

The organization of the paper is as follows. In section
\ref{sec:dynqgroups} we recall the algebraic definitions for
dynamical quantum groups such as \halg s,
\hHopfalg s, etc, and we introduce notation. Moreover, we recall
the FRST-construction in this setting and we consider the example
of the elliptic $U(2)$ dynamical quantum group. In section
\ref{sec:pairing} we recall the definition of a pairing between two
\hHopfalg s, and we extend this to include a $\ast$-structure and
we work out the details for the main example. In
section \ref{sec:actionsfrompairings} we introduce actions of one
\hHopfalg\ on another one in case there exists a pairing. In
section \ref{sec:repres} we apply this to representations and
corepresentations, and work out the details for
the main example. In particular, we calculate the pairing of two
matrix elements of irreducible corepresentations in terms
of elliptic hypergeometric series.
In section \ref{sec:singular} we propose
definitions of singular and spherical vectors, which we
intend to apply in more elaborate examples such as analogues of
dynamical symmetric spaces in future work.

\noindent
{\bf Acknowledgement.} We thank Hjalmar Rosengren for
useful discussions, and in particular for sharing
some unpublished results. We thank Peter Schauenburg for
sending \cite{Scha}.

\section{Preliminaries on dynamical quantum groups}
\label{sec:dynqgroups}

The algebraic notion for dynamical quantum groups is
in terms of \hHopfalg s. We recall these notions in
\S \ref{ssec:algebraicnotions}. Note that in case
the vector space $\h$ involved equals the trivial
space $\{0\}$, we are back in the case of Hopf algebras.
In \S \ref{ssec:FRST-construction} we recall the basic
FRST-construction for \hbialg s, and we discuss the
\hHopfalg\ extension for the explicit case of the
elliptic $U(2)$ dynamical quantum group.
In \S \ref{sec:singular} we also consider
$\h$-subalgebroids etc. for different vector spaces $\h$.

\subsection{Algebraic notions}\label{ssec:algebraicnotions}

Dynamical quantum groups have been introduced by
Etingof and Var\-chen\-ko \cite{EtinV}, see also the
lecture notes by Etingof and Schiffmann \cite{EtinS}.
For \hHopfalg s we follow the slightly different
definition as in \cite{KoelR}, see also \cite{Rose},
and we discuss shortly other approaches at the end of this
subsection in Remark \ref{rmk:comparison}.

We denote by $\hs$ a finite-dimensional complex vector
space. The notation is influenced by a natural construction
where $\h$ occurs as the Cartan subalgebra of a semisimple
Lie algebra, see \cite{EtinS}. By $\Mh$ we denote the
field of meromorphic functions on $\hs$, and the
function identically equal to $1$ on $\hs$ is denoted
by $\1 \in \Mh$.

A {\em \hprealg}\ $\A$ is a complex vector space
equipped with a decomposition $\A=\bigoplus_{\al,\be\in\hs}
\A_{\al,\be}$ and two left actions
$\mu_l^\A,\mu_r^\A \colon \Mh\to \text{End}_\C(\A)$, the
{\em left and right moment map}, which preserve the decomposition
and such that $\mu_l^\A(f)\mu_r^\A(g)=\mu_r^\A(g)\mu_l^\A(f)$
for all $f,g\in\Mh$. A {\em \hprealg\ homomorphism} is a
$\C$-linear map preserving the moment maps and the
decomposition.

A {\em \halg}\ $\A$ is a \hprealg , which is
also a unital associative algebra such that
the decomposition is a
bigrading  for the algebra, i.e.
$m^\A\colon \A_{\al,\be}\times \A_{\ga,\de} \to
\A_{\al+\ga,\be+\de}$ where $m^\A$ denotes the multiplication
of $\A$, and such that the
left and right moment map
$\mu_l^\A,\mu_r^\A \colon \Mh\to\A_{00}$ given by
$\mu_l^\A(f) = \mu_l^\A(f)1_\A$, $\mu_r^\A(f) = \mu_r^\A(f)1_\A$
are algebra embeddings, and such that
the commutation relations
\begin{equation}\label{eq:deflrmomentmapcommut}
\mu_l^\A(f)\, a = a \, \mu_l^\A(T_\al f),
\quad \mu_r^\A(f)\, a = a \, \mu_r^\A(T_\be f),
\qquad \forall \, a\in\A_{\al,\be}, \ \forall \, f\in\Mh,
\end{equation}
hold, where $T_\al$ is the automorphism of $\Mh$ defined
by $(T_\al f)(\la)=f(\la+\al)$. Note that in case $\A=\A_{00}$
this is just an extension of scalars to $\Mh$.
A {\em \halg\ homomorphism}
$\phi\colon\A\to\B$ of \halg s is an algebra homomorphism
which preserves the moment maps and the bigrading, i.e.
$\phi(\mu_l^\A(f))= \mu_l^\B(f)$,
$\phi(\mu_l^\A(f))= \mu_l^\B(f)$ and
$\phi(\A_{\al,\be})\subseteq \B_{\al,\be}$.

An important example is the algebra $\Dh$ of finite
difference operators $\sum_i f_i T_{\al_i}$
on the space $\Mh$ of meromorphic functions on $\h^\ast$.
Then $\Dh$ is a \halg .
The grading is given by $f T_{-\al}\in (\Dh)_{\al\al}$,
so $(\Dh)_{\al\be}=\{0\}$ for $\al\not=\be$. The
left and right moment map are equal to the natural embedding
by viewing an element $f\in \Mh$ as a multiplication operator
by $f$.
For later use we also note the identity
\begin{equation}\label{eq:trividinDh}
RS\1 =  (R\1) (T_{-\al}S\1), \qquad R\in (\Dh)_{\al\al}, S\in\Dh.
\end{equation}

Note that for a \halg\ $\A$ we obtain the \halg\ $\A^\lr$ by
interchanging the left and right moment maps. So
$(\A^\lr)_{\al,\be} = \A_{\be,\al}$, $\mu_l^{\A^\lr}=\mu^\A_r$,
$\mu_r^{\A^\lr}=\mu^\A_l$. Also, the opposite algebra $\A^\opp$,
i.e. $m^{\A^\opp}= m^\A\circ P$ with $P\colon \A\times\A\to\A\times \A$,
$(a, b) \mapsto (b, a)$ the natural flip operator, is again a \halg\  with
$\mu_l^{\A^\opp}=\mu_l^\A$, $\mu_r^{\A^\opp}=\mu_r^\A$,
$(\A^\opp)_{\al,\be}=\A_{-\al,-\be}$.

The {\em (matrix) tensor product} $\A \wtt \B$
of two \halg s $\A$ and $\B$ is a \halg \ with the
following definitions of the bigrading, moment maps
and multiplication;
\begin{subequations}\label{eq:defAtildetensorB}
\begin{align}
\label{eq:gradingAtildetensorB}
&\bigl( \A \wtt \B\bigr)_{\al\be} =\bigoplus_\ga
( \A_{\al\ga}\otimes_\Mh \B_{\ga\be} ),  \\
\label{eq:momentAtildetensorB}
&\mu_l^{\A \wtt \B}(f) = \mu_l^\A(f)\otimes 1, \qquad
\mu_r^{\A \wtt \B}(f) = 1\otimes \mu_r^\B(f), \\
\label{eq:productAtildetensorB}
&(a\otimes b)(c\otimes d) = (ac)\otimes (bd),
\end{align}
\end{subequations}
where $\otimes_\Mh$ denotes the tensor product
modulo the relations
\begin{equation}\label{eq:defotimesMh}
\mu_r^\A(f)a\otimes b =  a\otimes \mu_l^\B(f)b,
\qquad a\in\A,\ b\in \B, \ f\in \Mh.
\end{equation}
It is straightforward to check that
$\A\wtt \B$ is a \halg . Also note that we can
define the (matrix) tensor product $\A\wtt\B$ of two
\hprealg s $\A$, $\B$ in the same way
but only requiring \eqref{eq:gradingAtildetensorB},
\eqref{eq:momentAtildetensorB}. It can be
checked that in this case $\A\wtt\B$ gives again
a \hprealg .

For later use we note that, as \hprealg s,
$\A \wtt \Dh \cong \A \cong \Dh \wtt \A$
by
$a = a \otimes T_{-\be} = T_{-\al}\otimes a$
for $a\in A_{\al\be}$
and using \eqref{eq:defotimesMh} this implies
\begin{equation}\label{eq:identDhotimesAwithA}
\mu_r^\A(f) a = a \otimes fT_{-\be},
\qquad \mu_l^\A(f)a = fT_{-\al}\otimes a, \qquad a\in A_{\al\be}.
\end{equation}
This identification holds in particular for \halg s, and
for $\A=\Dh$ this gives the identification
$\Dh\wtt\Dh\cong \Dh$ via $fT_{-\al}\otimes gT_{-\al} = (fg) T_{-\al}$.

A {\em \hcoalg}\ is a \hprealg\ $\A$ with two
\hprealg\ homomorphisms $\De^\A\colon \A\to\A\wtt\A$,
the {\em comultiplication}, $\ep^\A\colon \A\to\Dh$, the {\em counit},
satisfying the coassociativity condition
$(\De^\A\otimes\id)\circ\De^\A = (\id\otimes\De^\A)\circ\De^\A$
and the counit condition
$(\ep^\A\otimes\id)\circ\De^\A= \id =
(\id\otimes\ep^\A)\circ \De^\A$ using the
identification \eqref{eq:identDhotimesAwithA}.

We use Sweedler's notation for the comultiplication, i.e.
$\De^\A(a)=\sum_{(a)} a_{(1)}\otimes a_{(2)}$ where the
decomposition on the right hand side is with respect to the
bigrading for a homogeneous element, so $a\in\A_{\al\be}$,
$a_{(1)}\in\A_{\al\eta}$, $a_{(2)}\in\A_{\eta\be}$. Then the
condition for the counit can be rewritten as
$\sum_{(a)}\mu_l^\A(\ep^\A(a_{(1)})\1)a_{(2)}=a =
\sum_{(a)}\mu_r^\A(\ep^\A(a_{(2)})\1)a_{(1)}$.
The counit is compatible with $\A\wtt\Dh\cong\A\cong  \Dh\wtt\A$
by  $\ep(a\otimes fT_{-\be}) =\ep^\A(a)\otimes \ep^\Dh (fT_{-\be})$
and $\ep(fT_{-\al}\otimes a) = \ep^\Dh(fT_{-\al})\otimes
\ep^\A(a)$ for $a\in\A_{\al\be}$.
Note that the maps $\ep^\A\otimes \ep^\Dh
\colon \A\wtt\Dh \to \Dh\wtt\Dh$ and
$\ep^\Dh\otimes \ep^\A\colon \Dh\wtt\A \to \Dh\wtt\Dh$
are well-defined.
So the counit axiom can be rewritten as
$\sum_{(a)} \ep^\A(a_{(1)})\otimes a_{(2)}\cong a
\cong \sum_{(a)}a_{(1)}\otimes\ep^\A(a_{(2)})$.

A {\em \hbialg}\  is a \halg\ $\A$ which is also
a \hcoalg\ and such that
the comultiplication  $\De^\A$ and
the counit $\ep^\A$ are \halg\ homomorphisms.
A {\em \hbialg\ homomorphism} $\phi\colon\A\to\B$ is a \halg\
homomorphism preserving the comultiplication and counit, i.e.
$(\phi\otimes \phi) \circ \De^\A = \De^\B\circ\phi$ and
$\ep^\B\circ\phi=\ep^\A$. Note that $\phi\otimes \phi\colon \A\wtt
\A \to \B\wtt \B$ is a well-defined operator.

We can make $\Dh$ into a \hbialg\ by setting
the comultiplication
$\De^\Dh \colon \Dh \to \Dh\wtt\Dh \cong\Dh$
to be the canonical isomorphism, and $\ep^\Dh$ to
be the identity.

A {\em \hHopfalg}\  is a \hbialg\ $\A$ equipped with
a $\C$-linear map, the antipode, $S^\A \colon \A \to \A$
satisfying $S^\A(\mu^\A_r(f)a) = S^\A(a) \mu^\A_l(f)$,
$S^\A(a\mu^\A_l(f)) = \mu^\A_l(f) S^\A(a)$ for all $a\in\A$ and
$f\in\Mh$ and
\begin{equation}\label{eq:defS}
\begin{split}
&m^\A\circ (\id\otimes S^\A)\circ \De^\A(a) = \mu_l^\A
(\ep^\A(a)\1), \qquad a \in \A, \\
&m^\A\circ (S^\A\otimes \id)\circ \De^\A(a) = \mu_r^\A
(T_\al \ep^\A(a)\1)), \qquad a \in \A_{\al,\be}.
\end{split}
\end{equation}
This definition follows \cite[\S 2]{KoelR}, and it differs
slightly from \cite{EtinV}, see also \cite{EtinS}. It is
straightforward to check that $m^\A\circ (S^\A\otimes \id)$ in
\eqref{eq:defS} is well-defined on $\A\wtt\A$, and for $m^\A\circ
(\id\otimes S^\A)$ we note that for $a\in\A_{\al\ga}$,
$b\in\A_{\ga\be}$, $f\in\Mh$, we have $m^\A\circ (\id\otimes
S^\A)((\mu_r^\A(f)a)\otimes b) = \mu_r^\A(f)aS^\A(b) =
aS^\A(b)\mu_r^\A(f) = a(S^\A(\mu_l^\A(f)b)= m^\A\circ (\id\otimes
S^\A)(a\otimes (\mu_l^\A(f)b))$ using $aS^\A(b)\in
\A_{\al\ga}\A_{-\be,-\ga} \subseteq \A_{\al-\be,0}$ by (iii)
below. With this definition it follows that antipode satisfies (i)
$S^\A(ab)=S^\A(b)S^\A(a)$, (ii) $\De^\A\circ S^\A = P\circ
(S^\A\otimes S^\A)\circ \De^\A$, (iii)
$S^\A(\A_{\al,\be})\subseteq \A_{-\be,-\al}$, (iv)
$S^\A(\mu^\A_l(f))=\mu^\A_r(f)$, $S^\A(\mu^\A_r(f))=\mu^\A_l(f)$,
(v) $S^\A(1)=1$, $\ep^\A\circ S^\A = S^{\Dh}\circ \ep$ where $P$
denotes the flip $P(a\otimes b)=b\otimes a$, and $S^\Dh\colon
\Dh\to\Dh$, $fT_{\al}\mapsto T_{-\al}\circ f = T_{-\al}\circ
(fT_{\al})\circ T_{-\al}$, is an algebra anti-isomorphism. See
\cite[Prop.~2.2]{KoelR}.

With the antipode as defined in the previous paragraph
the \hbialg\ $\Dh$ is a \hHopfalg .

The antipode on a \hHopfalg\ is compatible with $\A \wtt \Dh
\cong \A \cong \Dh \wtt \A$ by $S(a\otimes
fT_{-\be})=S^{\Dh}(fT_{-\be})\otimes S^\A(a)$ and
$S(fT_{-\al}\otimes a)=S^\A(a) \otimes S^{\Dh}(fT_{-\al})$ using
the notation of \eqref{eq:identDhotimesAwithA}. Note that these
maps are well-defined.

Note that for a \hbialg\ $\A$ the opposite \halg\ $\A^\opp$ is
again a \hbialg\  with $\De^{\A^\opp}=\De^\A$,
$\ep^{\A^\opp}=S^{\Dh}\circ\ep^\A$. Similarly, we can construct
the co-opposite \hbialg\ $\A^\cop$ which has the same algebra
structure, but with interchanged moment maps
$\mu_l^{\A^\cop}=\mu_r^\A$, $\mu_r^{\A^\cop}=\mu_l^\A$, and
$(\A^\cop)_{\al,\be}=\A_{\be,\al}$. (So as an \halg\ we have
$\A^\cop=\A^\lr$.) Moreover, the comultiplication is defined by
$\De^{\A^\cop}=P\circ\De^\A$, and the counit
$\ep^{\A^\cop}=\ep^\A$. Furthermore, if $\A$ is a \hHopfalg\ with
invertible antipode, then $\A^\opp$ and $\A^\cop$ are \hHopfalg s
with $S^{\A^\opp}=(S^\A)^{-1}$ and $S^{\A^\cop}=(S^\A)^{-1}$.

Assume that $\hs$ is equipped with a conjugation
$\la\mapsto\bar\la$, or
equivalently, a real form. Then this defines an
operator on $\Mh$ by $\bar f(\la)=\overline{f(\bar\la)}$.
A $\ast$-operator on a \halg\ $\A$ is a $\C$-antilinear
antimultiplicative involutive map $a\mapsto a^\ast$, such that
$\mu_l^\A(f)^\ast=\mu_l^\A(\bar f)$,
$\mu_r^\A(f)^\ast=\mu_r^\A(\bar f)$. This implies
that $(\A_{\al\be})^\ast=\A_{-\bar\al,-\bar\be}$.
In this case we see that $\Dh$ has a $\ast$-structure
given by $(fT_\al)^\ast= (T_{-\bar\al}f)T_{-\bar\al} =
T_{-\bar\al}\circ \bar f$, viewing $f\in\Dh$ as
the operator of multiplication by $f$. Note that the antipode
and the $\ast$-operator in $\Dh$ commute, which is not true
for a general \hHopfster .
If $\A$ is a \hbialg\ we require
the $\ast$-structure  to satisfy
$(\ast\otimes\ast)\circ\De^\A= \De^\A\circ\ast$ and
$\ep^\A\circ\ast= \ast\circ\ep^\A$, with the
$\ast$-operator on the right hand side the $\ast$-operator
of $\Dh$.  If $\A$ is a
\hbialg\ with a $\ast$-structure, and if moreover $\A$
is a \hHopfalg\ with invertible antipode, then
$S\circ\ast$ is an involution \cite[Lemma~2.9]{KoelR}.
We call $\A$ a \hHopfster .
The $\ast$-structure on $\A$ is compatible with
$\A \wtt \Dh \cong \A \cong \Dh \wtt \A$ by
$(a\otimes fT_{-\be})^\ast =a^\ast\otimes (fT_{-\be})^\ast$ and
$(fT_{-\al}\otimes a)^\ast =  (fT_{-\al})^\ast\otimes a^\ast$.
Note that these maps are well-defined.

For $\A$ a \hHopfster\ with invertible antipode, $\A^\opp$ and
$\A^\cop$ are \hHopfster \ for the same $\ast$-structure.

\begin{remark}\label{rmk:comparison}
Related structures motivated by the quantisation of
Poisson-Lie groupoids  have appeared before in the papers
by Lu \cite{Lu} and  Xu \cite{Xu}, see
also the discussion and motivation
in Etingof and Schiffmann \cite{EtinS}.  It has been proved by
Brzezi\'nski and Militaru \cite{BrzeM} that on the
bialgebroid level these structures are equivalent to
Takeuchi's \cite{Take} notion of $\times_R$-bialgebras.
Since Etingof and Nikshych \cite{EtinN} have proved that
weak Hopf algebras fit into the framework of Lu \cite{Lu},
this also includes weak bialgebras.

We can also view a \hbialg\ as a $\times_R$-bialgebra
equipped with a grading. Indeed, using the notation as
in Schauenburg \cite[\S 2]{Scha}, we see that $\Mh$ and
$\Dh$ play the
roles of $R$ and $\text{End}(R)$, where the action
of $R$ is given by the left moment map and the
action of $\bar R = R^\opp(=R)$ is given by the
right moment map. Moreover,
the tensorproduct \eqref{eq:defotimesMh} corresponds
to the notation $\int_r {}_{\bar r}\A\otimes {}_r\B$
and the weight space requirement \eqref{eq:gradingAtildetensorB} in the
tensor product corresponds to
$\int^r \A_{\bar r} \otimes \B_{\bar r}$, so that
Takeuchi's notion $\A \times_R \B$ corresponds to the (matrix)
tensor product. However, the antipode in the context of
$\times_R$-bialgebras is much more involved, and it seems that
this has not been settled yet, see Schauenburg \cite[\S 3]{Scha}
for a discussion based on categorical concepts and
B\"ohm and Szlach\'anyi \cite{BohS} for a
somewhat more restrictive but algebraic definition.
See also the references in the papers \cite{BohS}, \cite{Scha}
for more information. Note that the antipode in
\hHopfalg\ is different, since it makes essential use
of the bigrading.
\end{remark}

\subsection{FRST-construction
and the elliptic  $U(2)$ dynamical quantum
group}\label{ssec:FRST-construction}

We recall the FRST-construction for the case of an $R$-matrix with
dynamical and spectral parameter, see
Etingof and Varchenko \cite{EtinV}.
The FRST-construction associates to a $R$-matrix a \hbialg ,
which in many cases can be extended to a \hHopfalg . In general,
the $R$-matrix does not need to satisfy the
quantum dynamical Yang-Baxter equation
\eqref{eq:QDYBE}, but for the most interesting examples this is
the case. We give the
details for the elliptic $R$-matrix for $\mathfrak{sl}(2)$,
and we refer to Felder and Varchenko \cite{FeldV} and \cite{KoelvNR}
for an even more explicit description of the corresponding
\hHopfalg .

We start by recalling the definition of the quantum dynamical
Yang-Baxter equation. Let $\h$ be a finite dimensional
complex vector space viewed as a commutative Lie algebra
as before, and
let $V=\bigoplus_{\al\in\hs} V_\al$ be a diagonalisable
$\h$-module. A meromorphic function $R\colon \hs\times \C \to
\text{End} (V\otimes V)$ is a $R$-matrix if
it is $\h$-invariant, i.e. commutes with the $\h$-action on
$V\otimes V$, and if it satisfies
the quantum dynamical Yang-Baxter equation (with spectral
parameter);
\begin{equation}\label{eq:QDYBE}
R^{12}(\la-h^{(3)}, z_{12})R^{13}(\la, z_{13})
R^{23}(\la-h^{(1)}, z_{23}) =R^{23}(\la,
z_{23})R^{13}(\la-h^{(2)}, z_{13}) R^{12}(\la, z_{12}).
\end{equation}
Here $z_{ij}=z_i/z_j$, and e.g. $R^{12}(\la-h^{(3)}, z_{12})
(u\otimes v\otimes w) =
\bigl( R(\la-\mu,z_{12})(u\otimes v)\Bigr)\otimes w$
where $w\in V_{\mu}$.

Let $V=\bigoplus_{\alpha\in\h^*}V_\alpha$ be a finite-dimensional
diagonalizable $\h$-module and $R\colon\hs \times\C \to
\mathrm{End}_\h(V\otimes V)$ a meromorphic function, so
$R(\la,z)$ commutes
with the $\h$-action on $V\otimes V$. Let $\{e_x\}_{x\in X}$ be a
homogeneous basis of $V$, where $X$ is an index set. Write
$R^{ab}_{x y}(\lambda,z)$ for the matrix elements of $R$,
\begin{eqnarray*}
R(\lambda,z)(e_a\otimes e_b)=\sum_{x,y\in X} R^{a b}_{x
y}(\lambda, z) e_x\otimes e_y,
\end{eqnarray*}
and define $\omega\colon X\to \hs$ by $e_x\in V_{\omega(x)}$.
Let $\A_R$ be the unital complex associative algebra generated by the
elements $\{L_{xy}(z)\}_{x,y\in X}$, with $z\in\C$, together with
two copies of $\Mh$, embedded as subalgebras. The elements of
these two copies are denoted by $\mu_l^{\A_R}(f)= f(\lambda)$
and $\mu_r^{\A_R}(f)=f(\mu)$,
respectively. The defining relations of $\A_R$ are
$f(\lambda)g(\mu) = g(\mu)f(\lambda)$,
\begin{equation}
\begin{split}
f(\lambda)L_{x y}(z) = L_{x y}(z) f(\lambda+\omega(x)),\quad
f(\mu)L_{x y}(z) =L_{x y}(z) f(\mu+\omega(y)),
\end{split}
\end{equation}
for all $f$, $g\in \Mh$, together with the $R L L$-relations
\begin{eqnarray}\label{eq:RLL}
\sum_{x,y\in X} R^{x y}_{a c}(\lambda, z_1/z_2)L_{x b}(z_1)L_{y
d}(z_2)= \sum_{x,y\in X} R_{x y}^{b d}(\mu, z_1/z_2)L_{c
y}(z_2)L_{a x}(z_1),
\end{eqnarray}
for all $z_1, z_2\in\C$ and $a, b, c, d\in X$.

The bigrading on $A_R$ is defined by $L_{x y}(z)\in
\A_{\omega(x),\omega(y)}$.
The $\h$-invariance of $R$ ensures
that the bigrading is compatible with the $RLL$-relations
\eqref{eq:RLL}.
The counit and comultiplication defined by
\begin{equation}\label{eq:epandDeFRST}
\ep^{\A_R}(L_{a b}(z))=\delta_{a b} T_{-\omega(a)}, \quad
\De^{\A_R}(L_{a b}(z))=\sum_{x\in X} L_{a x}(z)\otimes L_{x b}(z)
\end{equation}
make $\A_R$ into a \hbialg , see
\cite{EtinV}.

\begin{example}\label{ex:ellipticunitary}
Take $\h\cong\hs\cong\C$ and let $V$ the two-dimensional $\h$-module
$V=\C e_1\oplus \C e_{-1}$. In the basis $e_1\otimes e_1$,
$e_1\otimes e_{-1}$, $e_{-1}\otimes e_1$, $e_{-1}\otimes e_{-1}$
the $R$-matrix is given by
\begin{equation}\label{eq:ellR}
R(\la,z)=R(\la,z,p,q)=
\begin{pmatrix}
1 & 0 & 0 & 0 \\
 0 & a(\la,z) & b(\la,z) & 0 \\
0 & c(\la,z) & d(\la,z) & 0 \\
 0 & 0 & 0 & 1
\end{pmatrix},
\end{equation}
where we assume $0<q<1$ and
\begin{equation}\label{eq:defabcd}
\begin{split}
a(\la,z)=\frac{\theta(z,q^{2(\lambda+2)})}
{\theta(q^2z,q^{2(\lambda+1)})},&\qquad
 b(\la,z)=\frac{\theta(q^2, q^{-2(\la+1)}z)}
 {\theta(q^2z, q^{-2(\la+1)})},\\
 c(\la,z)=\frac{\theta(q^2,q^{2(\la+1)}z)}
 {\theta(q^2z,q^{2(\la+1)})},&\qquad
 d(\la,z)=\frac{\theta(z,q^{-2\la})}{\theta(q^2z,q^{-2(\la+1)})}.
\end{split}
\end{equation}
Here the theta functions are normalised theta functions defined
by
\begin{equation}\label{eq:thetafunction}
\theta(z) = (z,p/z;p)_\infty,\ \ z\in\C\backslash\{0\},
\qquad \theta(a_1,\ldots, a_r)
=\prod_{i=1}^r \theta(a_i),
\end{equation}
where we assume $0<p<1$ and we employ the notation
\begin{equation}\label{eq:infproducts}
(a;p)_k = \prod_{i=0}^{k-1} (1-ap^i), \qquad
(a;p)_\infty  = \lim_{k\to\infty} (a;p)_k, 
\qquad (a_1,\ldots, a_r;p)_k = \prod_{i=1}^r (a_i;p)_k.
\end{equation}
For later use we note that theta functions satisfy
$\theta(pz)=\theta(z^{-1})=-z^{-1}\theta(z)$ and
the following addition formula
\begin{equation}\label{eq:addell}
\theta(xy, x/y, zw, z/w)= \theta(xw,x/w,zy,z/y) +
(z/y)\theta(xz,x/z,yw,y/w).
\end{equation}

The $R$-matrix defined by \eqref{eq:ellR} satisfies the
quantum dynamical Yang-Baxter equation
\eqref{eq:QDYBE}, see e.g. \cite{FeldV}, \cite{JimbOKS}, \cite{KoelvNR},
and references given there.

The four $L$-generators are denoted by
$\al(z)=L_{1,1}(z)$, $\be(z)=L_{1,-1}(z)$,
$\ga(z)=L_{-1,1}(z)$ and $\de(z)=L_{-1,-1}(z)$.
We do not give the relations for the generators arising from
\eqref{eq:RLL} explicitly, but we refer to \cite{FeldV},
\cite{KoelvNR}.
The corresponding \hbialg\ contains a group-like central element,
$\det(z)=\mu_r(F)\mu_l(F^{-1})\left[ \al(z)\de(q^2
z)-\ga(z)\be(q^2z)\right]$, with
$F(\la)=q^{\la} \theta(q^{-2(\la+1)})$, and adjoining the
inverse, denoted by $\det^{-1}(z)$, gives a \hHopfalg\
structure with antipode given by
\begin{equation}
\begin{split}\label{eq:antipode}
S(\alpha(z))= \mu_r(F)\mu_l(F^{-1})
\det^{-1}(q^{-2}z)\de(q^{-2}z), \quad
& S(\be(z)) = -\mu_r(F)\mu_l(F^{-1})
\det^{-1}(q^{-2}z)\be(q^{-2}z),\\
S(\ga(z)) = - \mu_r(F)\mu_l(F^{-1})
\det^{-1}(q^{-2}z)\ga(q^{-2}z), \quad
& S(\de(z)) = \mu_r(F)\mu_l(F^{-1})
\det^{-1}(q^{-2}z)\al(q^{-2}z),\\
S(\det^{-1}(z))=\det(z). \quad&
\end{split}
\end{equation}

We denote the corresponding \hHopfalg\ by
$\E=\bigoplus_{k,l\in\Z} \E_{kl}$. It can be made
into \hHopfster\ by defining
$\det^{-1}(z)^\ast=\det^{-1}(q^{-2}/\bar z)$,
$\al(z)^\ast=\de(1/\bar z)$, $\be(z)^\ast=-\ga(1/\bar z)$,
$\ga(z)^\ast=-\be(1/\bar z)$, $\de(z)^\ast=\al(1/\bar z)$,
see \cite{KoelvNR}. We
come back to this example throughout this paper.
\end{example}

\begin{example}\label{ex:rationalgl2}
With the same convention for $\h$, $\hs$, $V$, $V\otimes V$ as
in Example \ref{ex:ellipticunitary} the $R$-matrix
\begin{equation}\label{eq:ratR}
R(\la)=R(\la,q)=
\begin{pmatrix}
q & 0 & 0 & 0 \\
0 & 1 & \frac{q^{-1}-q}{q^{2(\la+1)}-1} & 0 \\
0 & \frac{q^{-1}-q}{q^{-2(\la+1)}-1} &
\frac{(q^{2(\la+1)}-q^2)(q^{2(\la+1)}-q^{-2})}{(q^{2(\la+1)}-1)^2} & 0 \\
0 & 0 & 0 & 1
\end{pmatrix}.
\end{equation}
satisfies the quantum dynamical Yang-Baxter equation
\eqref{eq:QDYBE} and is also  $\h$-invariant.
The corresponding \halg\ arising from the
FRST-construction can also be made into a \hHopfalg .
This example can be found in \cite[\S 2.2]{EtinS}, and the
corresponding harmonic analysis has been studied in \cite{KoelR}.
It has been proved by Stokman \cite{Stok} that this case can
be obtained from the usual quantum $SL(2)$ group by use of a vertex
IRF transform based on the twisted coboundary element introduced
by Babelon, Bernard and Billey \cite{BabeBB}, see also
\cite{BuffR}, \cite{RoseCM} for more information.

The $R$-matrix in \eqref{eq:ratR} can be obtained by
a suitable limit transition from the $R$-matrix \eqref{eq:ellR}.
To see this we recall that the $R$-matrix in \eqref{eq:ellR}
is gauge equivalent to the elliptic $R$-matrix in \cite[\S 3]{JimbOKS} by
an explicit closed multiplicative $2$-form.
In the notation of \cite[\S\S 1-2]{EtinV}, \cite[\S 6.2]{EtinS}
we only have to specify the multiplicative $2$-form
$\displaystyle{\phi_{12}(\la)= \frac{1}{q} \frac{(pq^{2\la},q^{-2(\la+1)};p)_\infty}
{(q^{-2(\la+2)},pq^{2(\la+1)};p)_\infty}}$.
The closedness follows from $\displaystyle{\phi_{12}(\la)=
\frac{f_2(\la) f_1(\la+1)}{f_1(\la)f_2(\la-1)}}$ with $f_1(\la)=1$,
$f_2(\la)=  q^{-\la} \bigl( q^{-2(\la+2)},pq^{2(\la+1)};p)_\infty\bigr)^{-1}$.
The limit transition from the $R$-matrix of \cite{JimbOKS}
by taking $\lim_{p\to 0}$, $\lim_{z\to 0}$ gives $q$ times the
$R$-matrix in \eqref{eq:ratR} with $q$ replaced by $q^{-1}$.
Then the map
$\displaystyle{L_{ab}(z) \mapsto \frac{\mu_r(f_b)}{\mu_l(f_a)} \tilde L_{ab}(z)}$
defines a \hbialg -isomorphism of the \hbialg\ from the FRST-construction
for the $R$-matrix \eqref{eq:ellR} to the \hbialg\ from the FRST-construction
for the $R$-matrix in \cite[\S 3]{JimbOKS}.
Now the limit transition of the elliptic $R$-matrix in
\cite[\S 3]{JimbOKS} to \eqref{eq:ratR} gives a corresponding
formal limit transition at the level of
the corresponding \hbialg s  and \hHopfalg s.
The harmonic analyis on the
corresponding \hbialg\ and \hHopfalg\ is studied in
\cite{KoelR}, and the pairing for this case is studied extensively
by Rosengren \cite{Rose}.
\end{example}

\section{Pairings for dynamical quantum groups}\label{sec:pairing}

We consider in this section a natural pairing between
dynamical quantum group. For certain
weak Hopf algebras this has been considered by Etingof and
Nikshych \cite[\S 5]{EtinN}. The pairing
is studied by Rosengren \cite{Rose} in the context
of duals of \hbialg s and \hHopfalg s. Rosengren \cite{Rose}
shows in particular that the quantized universal
enveloping algebra $U_q(\mathfrak{sl}(2))$ is contained
in the \hHopfalg\ associated to the $R$-matrix \eqref{eq:ratR}.
We study the pairing for \hHopfster s and recall the
explicit pairings for the \hbialg s arising from the FRST-construction.

The following definition has been introduced by
Rosengren \cite[\S 3]{Rose}. This definition is
different from the pairing for $\times_R$-bialgebras
as introduced by Schauenburg \cite[\S 5]{Scha}, cf.
Remark \ref{rmk:comparison}.

\begin{defn}\label{def:pairingforhbialgebras}
A pairing for \hbialg s $\U$ and $\A$ is a $\C$-bilinear map
$\langle \cdot,\cdot\rangle \colon \U\times \A \to \Dh$ satisfying
\begin{subequations}\label{eq:pairingforhbialgebras}
\begin{align}
\label{eq:pairing-target}
& \pair{\U_{\al\be}}{\A_{\ga\de}}
\subseteq (\Dh)_{\al+\de,\be+\ga}, \\
\label{eq:pairing-leftmoment}
& \pair{\mu_l^\U(f)X}{a} = \pair{X}{\mu_l^\A(f) a}
 = f \circ \pair{X}{a}, \\
 \label{eq:pairing-rightmoment}
& \pair{X\mu_r^\U(f)}{a}= \pair{X}{a\mu_r^\A(f)}
 = \pair{X}{a} \circ f, \\
\label{eq:pairing-leftproduct}
&\pair{XY}{a}=\sum_{(a)} \pair{X}{a_{(1)}}T_\rho \pair{Y}{a_{(2)}},
    \quad \De^\A(a)=\sum_{(a)} a_{(1)}\otimes a_{(2)},\;
    a_{(1)}\in \A_{\ga\rho},\\
\label{eq:pairing-rightproduct}
&\pair{X}{ab}=\sum_{(X)} \pair{X_{(1)}}{a}T_\rho\pair{X_{(2)}}{b},
     \quad \De^\U(X)=\sum_{(X)} X_{(1)}\otimes X_{(2)},\;
     X_{(1)}\in \U_{\al\rho},\\
\label{eq:pairing-counit-unit}
&\pair{X}{1}=\ep^\U(X), \qquad \pair{1}{a}=\ep^\A(a).
\end{align}
\end{subequations}
If moreover, $\U$ and $\A$ are \hHopfalg s, then we
require
\begin{equation}\label{eq:pairingforhHopfalgebras}
\pair{S^\U(X)}{a} = S^\Dh( \pair{X}{S^\A(a)}).
\end{equation}
\end{defn}

\begin{remark}\label{rmk:pairingforhbialgebras}
(i) Note \eqref{eq:pairing-target} implies that
$\pair{X}{a}=0$ whenever $X\in\U_{\al\be}$, $a\in\A_{\ga\de}$
with $\al+\de\not=\be+\ga$. \par\noindent
(ii) In \eqref{eq:pairing-leftmoment} and
\eqref{eq:pairing-rightmoment} we consider $f$ as
a multiplication operator in $\Dh$. \par\noindent
(iii) Requiring only \eqref{eq:pairing-target},
\eqref{eq:pairing-leftmoment}, and \eqref{eq:pairing-rightmoment}
gives the definition of a pairing between \halg s $\U$ and $\A$.
Similarly, we can define a pairing between a \halg \ $\U$
and \hcoalg\ $\A$ by requiring
\eqref{eq:pairing-target}, \eqref{eq:pairing-leftmoment},
\eqref{eq:pairing-rightmoment}, with the convention
that in the \hcoalg\ $a\mu_r^\A(f) = \mu_r(T_{-\de}f)a$
for $a\in\A_{\ga\de}$, \eqref{eq:pairing-leftproduct}
and the last equality of \eqref{eq:pairing-counit-unit}.
\par\noindent
(iii) Note that \eqref{eq:pairingforhHopfalgebras} is
self-dual, since the antipode of $\Dh$ is involutive.
If a pairing exists for $\U$ and $\A$ and assuming
the pairing is non-degenerate, and $\A$ is \hHopfalg\
we can use \eqref{eq:pairingforhHopfalgebras} to
equip $\U$ with an antipode. In the case $\U$ and $\A$
are \hHopfalg s, which are \hbialg s with a pairing,
it follows from \eqref{eq:pairing-leftproduct},
\eqref{eq:pairing-rightproduct} and the antipode being
antimultiplicative that it suffices to check
\eqref{eq:pairingforhHopfalgebras} for generators.
\end{remark}

For applications to dynamical quantum groups we are
particularly interested in a specific type of pairings
induced by $R$-matrices as studied previously
in \cite{EtinN}, \cite{Rose}. For this we need the
following definition, see \cite[Def.~3.16]{Rose}, and
recall the definition of the coopposite \hbialg\
as in \S \ref{ssec:algebraicnotions}.

\begin{defn}\label{defn:cobraiding}
A cobraiding on a \hbialg\ $\A$ is a
pairing $\pair{\cdot}{\cdot}\colon\A^\cop \times
\A\to\Dh$ for \hbialg s satisfying
\begin{equation}\label{eq:defcobraiding}
\sum_{(a),(b)} \mu_l^\A(\pair{a_{(1)}}{b_{(1)}}\1)
a_{(2)}b_{(2)} =
\sum_{(a),(b)} \mu_r^\A(\pair{a_{(2)}}{b_{(2)}}\1)
b_{(1)}a_{(1)},
\end{equation}
as an identity in $\A$ and where
$\De^\A(a)=\sum_{(a)} a_{(1)}\otimes a_{(2)}$,
$\De^\A(b)=\sum_{(b)} b_{(1)}\otimes b_{(2)}$.
\end{defn}

\begin{remark}\label{rmk:cobraidingandRmatrix}
{(i)} It suffices to check \eqref{eq:defcobraiding} for generators
of the \hbialg\ $\A$, see \cite[Lemma 3.17]{Rose}. For this reason we
need the coopposite \halg\  in the first leg of the pairing.
\par\noindent
{(ii)} The relation \eqref{eq:defcobraiding} for $a$ and
$b$ generators of a \hbialg\ constructed
by the FRST-construction as in \S \ref{ssec:FRST-construction}
can be matched with
the quadratic relations \eqref{eq:RLL} using \eqref{eq:epandDeFRST}.
Then the fact that this
pairing defined on the generators in terms of the
$R$-matrix by
\begin{equation*}
\pair{L_{ij}(w)}{L_{kl}(z)} = R^{jl}_{ik}(\la,\frac{w}{z})
\, T_{-\om(i)-\om(k)}
\end{equation*}
extends to a pairing on $(\A_R)^\cop\times \A_R$ satisfying the
conditions of Definition \ref{def:pairingforhbialgebras} is
equivalent to the $R$-matrix being a solution to the quantum
dynamical Yang-Baxter equation \eqref{eq:QDYBE}.
Consequently, any solution of
the quantum dynamical Yang-Baxter equation gives rise to a
\hbialg\ with a cobraiding. See Rosengren \cite[\S 3]{Rose} for
these results in case there is no spectral parameter, which can
easily be adapted to include the case of an $R$-matrix with
both dynamical parameters and a spectral parameter.
\par\noindent
{(iii)} In general the cobraiding arising from an $R$-matrix
is not non-degenerate. In case of a finite dimensional weak Hopf algebra
the pairing is non-degenerate, see Etingof and Nikshych \cite{EtinN}.
For the case of Example \ref{ex:rationalgl2} the radical is non-trivial,
and it has been determined explicitly by Rosengren \cite[\S 5]{Rose}.
\end{remark}

For later use we also need the relation between $\ast$-structures
on \hHopfalg s and pairings. It turns out that we need
this in two, closely related, variants, one of them well-suited
for left corepresentations and the other well-suited for
right corepresentations. We assume that one of
the \hHopfster s has a fixed $\ast$-operator, and for the
other we use the notation $\ast$ and $\da$
depending on the pairing. This definition has been
communicated to us by Hjalmar Rosengren.

\begin{defn}\label{def:pairingster} Let $\U$, $\A$ be \hHopfalg s
with invertible antipodes  being paired as \hHopfalg s.
Assume that $\A$ is a \hHopfster . We say that
$\U$ and $\A$ are paired as \hHopfster s if moreover
$\U$ is a \hHopfster\ and
\begin{equation}\label{eq:pairingXstara}
\pair{X^\ast}{a} = T_{-\ga}\circ
(\pair{X}{S^\A(a)^\ast})^\ast\circ T_{-\de}, \qquad
\forall\, a\in\A_{\ga\de}, \ \forall \, X\in\U,
\end{equation}
or if $\U$ is a \hHopfster\ (with $\ast$ denoted by $\da$) and
\begin{equation}\label{eq:pairingXdaggera}
\pair{X^\da}{a} = T_{-\ga}\circ
(\pair{X}{S^\A(a^\ast)})^\ast\circ T_{-\de}, \qquad
\forall\, a\in\A_{\ga\de}, \ \forall \, X\in\U .
\end{equation}
\end{defn}

It will be clear from the context and/or the notation
in the sequel if we use \eqref{eq:pairingXstara}
or \eqref{eq:pairingXdaggera} for paired \hHopfster s.
For pairings for Hopf algebras the standard choice is
compatible with \eqref{eq:pairingXstara}.

\begin{remark}\label{rmk:defpairingster}
(i) It is not clear from \eqref{eq:pairingXstara} and
\eqref{eq:pairingXdaggera} that
$\ast$ and $\da$ on $\U$ are compatible with the
usual properties of a $\ast$-operator. This is
the content of Proposition \ref{prop:astandpairing}.
\par\noindent
(ii) If $\U$ and $\A$ are paired as \hHopfalg s with
invertible antipodes, and $\A$ is also a \hHopfster , it
follows that $\pair{(X^\da)^\ast}{a}=
\pair{X}{S^2(a)}= \pair{S^2(X)}{a}$, or the
\hHopfalg\ isomorphism $S^2\colon\U\to\U$ interwines
the two $\ast$-structures, $S^2(X^\da) = X^\ast$.
Note that this implies $S\circ\da = \ast\circ S$.
The choice for the $\ast$-structure on $\U$
is related to the notion
of unitarisability for right or left
corepresentations, see \S \ref{sec:repres}.
\par\noindent
(iii) Assume that $\U$ and $\A$ are \hHopfster s
with invertible antipodes that are paired
as in \eqref{eq:pairingXstara} or
\eqref{eq:pairingXdaggera}. Note that this
definition is not obviously self-dual, but it is:
\begin{equation}\label{eq:selfpairingdefpairing}
\begin{split}
\pair{X}{a^\ast} &=  T_{-\al}\circ \bigl(
\pair{S^\U(X)^\ast}{a}\bigr)^\ast \circ T_{-\be}, \\
\pair{X}{a^\ast} &=  T_{-\al}\circ \bigl(
\pair{S^\U(X^\da)}{a}\bigr)^\ast \circ T_{-\be},
\end{split}
\end{equation}
for $X\in\U_{\al\be}$, $a\in\A$ and assuming the
pairing \eqref{eq:pairingXstara} for the
first equation of \eqref{eq:selfpairingdefpairing}
and the
pairing \eqref{eq:pairingXdaggera} for the
second equation of \eqref{eq:selfpairingdefpairing}.
We prove the first equation of \eqref{eq:selfpairingdefpairing},
the second being proved more easily.
The assumption
that $\A$ and $\U$ have invertible antipodes implies
$S\circ \ast$ is an involution.
For $X\in\U_{\al\be}$, $a\in\A_{\ga\de}$ we have
$a^\ast\in\A_{-\bar\ga,-\bar\de}$, so that
\begin{equation*}
\begin{split}
\pair{X}{a^\ast} &= \pair{(X^\ast)^\ast}{a^\ast}
= T_{\bar\ga} \circ (\pair{X^\ast}{S^\A(a^\ast)^\ast})^\ast \circ T_{\bar \de}
= T_{\bar\ga} \circ (\pair{X^\ast}{(S^\A)^{-1}(a)})^\ast \circ T_{\bar \de}
\\
&= T_{\bar\ga} \circ \bigl( S^\Dh(\pair{(S^\U)^{-1}(X^\ast)}{a})\bigr)^\ast
\circ T_{\bar \de}
=T_{\bar\ga} \circ \bigl( S^\Dh(\pair{S^\U(X)^\ast}{a})\bigr)^\ast
\circ T_{\bar \de}
\end{split}
\end{equation*}
using \eqref{eq:pairingforhHopfalgebras}
and $S^\Dh$ being an involution as well. Write
$\pair{S^\U(X)^\ast}{a}=fT_{-\bar\be-\de}=fT_{-\bar\al-\ga}$
with $f=\pair{S^\U(X)^\ast}{a}\1$ for the moment, so
\begin{equation*}
\pair{X}{a^\ast} 
= T_{\bar\ga} \circ \bar f \circ
T_{-\be-\bar\de}\circ T_{\bar\de} =
T_{-\al} \circ (fT_{-\bar\al-\ga})^\ast\circ T_{-\be} 
= T_{-\al} \circ \bigl( \pair{S^\U(X)^\ast}{a} \bigr)^\ast\circ
T_{-\be},
\end{equation*}
proving the first equality of \eqref{eq:selfpairingdefpairing}.
\end{remark}

\begin{prop}\label{prop:astandpairing}
Assume that $\U$ and $\A$ are paired \hbialg s, and
that $\A$ is a \hHopfster\ with invertible antipode.
Define the pairing $\pair{X^\ast}{a}$ by \eqref{eq:pairingXstara},
then for all $X,Y\in\U$, for all
$a,b\in\A$, for all $c_1,c_2\in\C$,
\begin{subequations}\label{eq:prop-astandpairing}
\begin{align}
\label{eq:Clinearity-prop:astandpairing}
&\pair{(c_1X+c_2Y)^\ast}{a} =
\overline{c_1} \pair{X^\ast}{a} +
\overline{c_2} \pair{Y^\ast}{a}, \\
\label{eq:antimultiplicativity-prop:astandpairing}
&\pair{(XY)^\ast}{a} = \pair{Y^\ast X^\ast}{a}, \\
\label{eq:involution-prop:astandpairing}
&\pair{(X^\ast)^\ast}{a} = \pair{X}{a}, \\
\label{eq:leftrightmoment-prop:astandpairing}
&\pair{(\mu^\U_l(f))^\ast}{a} = \pair{\mu_l^\A(\bar f)}{a},
\qquad
\pair{(\mu^\A_r(f))^\ast}{a} = \pair{\mu_r^\A(\bar f)}{a}, \\
\label{eq:counit-prop:astandpairing}
&\pair{X^\ast}{1} = (\ep^\U(X))^\ast, \\
\label{eq:comultiplication-prop:astandpairing}
&\pair{X^\ast}{ab} = \sum_{(X)} \pair{(X_{(1)})^\ast}{a} \circ
T_{-\bar\eta} \circ \pair{(X_{(2)})^\ast}{b}, \quad
\De^\U(X)=\sum_{(X)} X_{(1)}\otimes
X_{(2)}, \ X_{(1)}\in \U_{\al\eta},
\end{align}
\end{subequations}
where we assume that $X\in\U_{\al\be}$ in
\eqref{eq:comultiplication-prop:astandpairing}.
The results in \eqref{eq:prop-astandpairing} remain
valid upon replacing $\ast$ by $\da$ in the left hand
sides of all pairings and using \eqref{eq:pairingXdaggera}.
\end{prop}

Comparing with Definition
\ref{def:pairingforhbialgebras}, we see that we can
look upon \eqref{eq:Clinearity-prop:astandpairing},
\eqref{eq:antimultiplicativity-prop:astandpairing}
and \eqref{eq:involution-prop:astandpairing} as the
weak formulation of the $\ast$-operator being
a $\C$-antilinear antimultiplicative involution.
Then \eqref{eq:leftrightmoment-prop:astandpairing}
is the weak formulation of being a $\ast$-operator
on a \halg . Finally, \eqref{eq:counit-prop:astandpairing}
is the weak formulation of $\ep\circ\ast = \ast\circ\ep$, and
\eqref{eq:comultiplication-prop:astandpairing} is the
weak formulation of $(\ast\otimes\ast)\circ\De =\De\circ\ast$.

\begin{proof} The $\C$-linearity of the pairing and
\eqref{eq:pairingXstara} imply \eqref{eq:Clinearity-prop:astandpairing}.
To prove \eqref{eq:antimultiplicativity-prop:astandpairing}
assume $X\in \U_{\al\be}$, $Y\in\U_{\si\tau}$ so that
$XY\in \U_{\al+\si,\be+\tau}$ and $a\in\A_{\ga\de}$, so
$\De^\A(a)=\sum_{(a)} a_{(1)}\otimes a_{(2)}$, $a_{(1)}\in \A_{\ga\eta}$,
$a_{(2)}\in \A_{\eta\de}$.
Using $(\ast\otimes \ast)\circ \De = \De\circ \ast$,
$P\circ(S\otimes S)\circ \De = \De\circ S$ in a \hHopfalg\ and
\eqref{eq:pairing-leftproduct} we find
\begin{equation*}
\begin{split}
\pair{(XY)^\ast}{a} &= T_{-\ga}\circ
(\pair{XY}{S^\A(a)^\ast})^\ast\circ T_{-\de} \\
&= T_{-\ga}\circ \Bigl( \sum_{(a)} \pair{X}{S^\A(a_{(2)})^\ast} T_{\bar\eta}
\pair{Y}{S^\A(a_{(1)})^\ast} \Bigr)^\ast\circ T_{-\de} \\
&=\sum_{(a)}T_{-\ga}\circ (\pair{Y}{S^\A(a_{(1)})^\ast})^\ast\circ T_{-\eta}\circ
(\pair{X}{S^\A(a_{(2)})^\ast})^\ast \circ T_{-\de} \\
&= \sum_{(a)} \pair{Y^\ast}{a_{(1)}}\circ T_{\eta}\circ
\pair{X^\ast}{a_{(2)}} = \pair{Y^\ast X^\ast}{a},
\end{split}
\end{equation*}
since $S^\A(a_{(2)})^\ast \in \A_{\bar\de,\bar\eta}$.
The proof of \eqref{eq:involution-prop:astandpairing} is an
immediate consequence of the fact that $S^\A\circ \ast$ is an
involution in a \hHopfalg\ with invertible antipode, see
\S \ref{ssec:algebraicnotions}, and the definition of the
$\ast$-operator in $\Dh$. For \eqref{eq:leftrightmoment-prop:astandpairing}
we use for $a\in \A_{\ga\de}$
\begin{equation*}
\begin{split}
\pair{(\mu^\U_l(f))^\ast}{a} &=
T_{-\ga}\circ
(\pair{\mu_l^\U(f)}{S^\A(a)^\ast})^\ast\circ T_{-\de}
=T_{-\ga}\circ ( f\circ \ep^\A (S^\A(a)^\ast))^\ast\circ T_{-\de}\\
&=T_{-\ga}\circ S^\Dh(\ep^\A (a))\circ \bar f \circ T_{-\de}
= g \circ \bar f \circ T_{-\de},
\end{split}
\end{equation*}
where we write $\ep^\A(a)=gT_{-\ga}$. On the other hand
$\pair{\mu_l^\U(\bar f)}{a} = \bar f \circ \ep^\A(a) = \bar f \circ gT_{-\ga}$.
Since $\ep^\A(a)\in \bigl(\Dh\bigr)_{\ga\de}$, both expressions
are zero for $\ga\not=\de$, see Remark \ref{rmk:pairingforhbialgebras}(i),
and that for $\ga=\de$ they are equal. The
statement for the right moment map in \eqref{eq:leftrightmoment-prop:astandpairing}
is proved similarly.

Finally, to prove \eqref{eq:comultiplication-prop:astandpairing} we
take  $X\in \U_{\al\be}$, so $\De^\U(X) =\sum_{(X)}X_{(1)}\otimes X_{(2)}$,
with $X_{(1)}\in \U_{\al\eta}$, $X_{(2)}\in \U_{\eta\be}$, and $a\in \A_{\ga\de}$,
$b\in \A_{\si\tau}$, so that the left hand side of  \eqref{eq:comultiplication-prop:astandpairing} can be rewritten as
\begin{equation*}
\begin{split}
\pair{X^\ast}{ab} &= T_{-\ga-\si}\circ (\pair{X}{S^\A(ab)^\ast})^\ast \circ T_{-\de-\tau}
=T_{-\ga-\si}\circ (\pair{X}{S^\A(a)^\ast S^\A(b)^\ast})^\ast \circ T_{-\de-\tau}\\
&=T_{-\ga-\si}\circ \Bigl( \sum_{(X)} \pair{X_{(1)}}{S^\A(a)^\ast}\circ T_\eta\circ
\pair{X_{(2)}}{S^\A(b)^\ast} \Bigr)^\ast\circ T_{-\de-\tau} \\
&=\sum_{(X)}T_{-\ga-\si}\circ (\pair{X_{(2)}}{S^\A(b)^\ast})^\ast\circ T_{-\bar\eta} \circ (\pair{X_{(1)}}{S^\A(a)^\ast})^\ast \circ T_{-\de-\tau}\\
&=\sum_{(X)}T_{-\ga} \circ \pair{(X_{(2)})^\ast}{b}\circ T_{\tau}\circ T_{-\bar\eta}
\circ T_{\ga}\circ \pair{(X_{(1)})^\ast}{a} \circ T_{-\tau}.
\end{split}
\end{equation*}
Write
$\pair{(X_{(1)})^\ast}{a} = fT_{\bar \eta-\ga}\in \bigl( \Dh\bigr)_{\de-\bar\al,\ga-\bar\eta}$, and
$\pair{(X_{(2)})^\ast}{b} = gT_{\bar \eta-\tau}\in \bigl( \Dh\bigr)_{\tau-\bar\eta,\si-\bar\be}$, so that we can rewrite the summand
as
\begin{equation*}
T_{-\ga} \circ (g T_{\bar\eta-\tau}) \circ T_{\tau}\circ T_{-\bar\eta}
\circ T_{\ga}\circ (fT_{\bar\eta-\ga}) \circ T_{-\tau} = f \circ (T_{-\ga} g) \circ T_{\bar\eta-\ga-\tau}.
\end{equation*}
Since this is also equal to
\begin{equation*}
\pair{(X_{(1)})^\ast}{a}\circ
T_{-\bar\eta} \circ \pair{(X_{(2)})^\ast}{b} = f \circ T_{-\ga} \circ gT_{\bar\eta-\tau},
\end{equation*}
the last result follows.

The proofs for the $\ast$-structure $\dagger$ follow analogously.
\end{proof}

\subsection{Elliptic $U(2)$ dynamical quantum group}
\label{ssec:ellpairing}

With the notation of \S \ref{ssec:FRST-construction}
for the explicit $R$-matrix
in \eqref{eq:ellR} of Example
\ref{ex:ellipticunitary} we have from
Remark \ref{rmk:cobraidingandRmatrix}(ii)
\begin{equation}\label{eq:elldualgenerator}
\begin{array}{lll}
\pair{\al(w)}{\al(z)} =T_{-2}, &
\pair{\al(w)}{\de(z)} =a(\la,w/z)T_0, &
\pair{\be(w)}{\ga(z)} =b(\la,w/z)T_0, \\[2mm]
\pair{\ga(w)}{\be(z)} =c(\la,w/z)T_0, &
\pair{\de(w)}{\al(z)} =d(\la,w/z)T_0, &
\pair{\de(w)}{\de(z)} =T_2,
\end{array}
\end{equation}
and all other pairings between generators being zero. In order to
calculate the pairing with the inverse determinant, we calculate
using \eqref{eq:elldualgenerator}, the definition of the elliptic
determinant in \S \ref{ssec:FRST-construction} and Definition
\ref{def:pairingforhbialgebras} directly
\begin{equation}\label{eq:ellpairingdet}
\begin{split}
&\pair{\left( \begin{matrix} \al(w) &\be(w) \\ \ga(w) &\de(w)
\end{matrix}\right)}{\det(z)}=
q\frac{\theta(q^{-2}w/z)}{\theta(w/z)}
\left( \begin{matrix} T_{-1} & 0 \\ 0 & T_1
\end{matrix}\right), \\
&\pair{\det(w)}{\left( \begin{matrix} \al(z) &\be(z)
\\ \ga(z) &\de(z) \end{matrix}\right)}=
q\frac{\theta(w/z)}{\theta(q^2w/z)}
\left( \begin{matrix} T_{-1} & 0 \\ 0 & T_1
\end{matrix}\right).
\end{split}
\end{equation}
Most of the pairings are straightforward to compute, and we show how to
calculate $\pair{\de(w)}{\det(z)}$, which is the most
involved. Recall from Example
\ref{ex:ellipticunitary} that
$\det(z) = \mu_r(F)\mu_l(F^{-1})\bigl(  \al(z)\de(q^2z)-\ga(z)\be(q^2z)\bigr)$,
with $F(\la)=q^{\la} \theta(q^{-2(\la+1)})$.
Since $\de(w)\in (\E^\cop)_{-1,-1}$, $\be(w)\in (\E^\cop)_{-1,1}$,
$\De^{\E^\cop}(\de(w))= \be(w)\otimes \ga(w) + \de(w)\otimes \de(w)$
we obtain
\begin{equation*}
\begin{split}
\pair{\de(w)}{\al(z)\de(q^2z)} &=
\pair{\be(w)}{\al(z)} T_1 \pair{\ga(w)}{\de(q^2z)} +
\pair{\de(w)}{\al(z)} T_{-1} \pair{\de(w)}{\de(q^2z)} = d(\la, \frac{w}{z}) T_1, \\
\pair{\de(w)}{\ga(z)\be(q^2z)} &=
\pair{\be(w)}{\ga(z)} T_1 \pair{\ga(w)}{\be(q^2z)} +
\pair{\de(w)}{\ga(z)} T_{-1} \pair{\de(w)}{\be(q^2z)} \\
&=
b(\la,\frac{w}{z}) T_1 c(\la,\frac{w}{q^2z}) =
b(\la,\frac{w}{z}) c(\la+1,\frac{w}{q^2z}) T_1,
\end{split}
\end{equation*}
so that, using the explicit expressions \eqref{eq:defabcd}
and elementary properties of the theta function,
\begin{equation*}
\begin{split}
\pair{\de(w)}{\det(z)} &= \frac{1}{F(\la)} \bigl( d(\la,\frac{w}{z})
- b(\la,\frac{w}{z})c(\la+1,\frac{w}{q^2z})\bigr) F(\la+1) \, T_1\\
&= \frac{ \theta(\frac{w}{z}, \frac{w}{z}, q^{-2\la}, q^{2(\la+2)})
-\theta(q^2,q^{-2(\la+1)}\frac{w}{z},q^2,q^{2(\la+1)}\frac{w}{z})}
{q\theta(q^2\frac{w}{z},q^{-2(\la+1)}, \frac{w}{z},q^{2(\la+1)})}\, T_1
\end{split}
\end{equation*}
and applying \eqref{eq:addell} (with $x^2$, $y^2$, $w^2$, $z^2$ replaced
by $q^{-2\la}w/z$, $q^{2\la}w/z$, $q^{-2(\la+2)}w/z$,
$q^{2(\la+2)}w/z$) we see that the numerator simplifies to
a product of $4$ theta functions partially cancelling the
theta functions in the numerator. This proves one of the
eight statements in \eqref{eq:ellpairingdet}. All other statements
are proved in this way, and only the pairings with the $\de(w)$,
$\de(z)$ require the use of \eqref{eq:addell}.

From \eqref{eq:ellpairingdet} we see how to extend the pairing
to $\E^\cop\times \E$;
\begin{equation*}
\left( \begin{matrix} T_{-1} & 0 \\ 0 & T_1
\end{matrix}\right) =
\pair{\left( \begin{matrix} \al(w) &\be(w) \\ \ga(w) &\de(w)
\end{matrix}\right)}{1} =
\pair{\left( \begin{matrix} \al(w) &\be(w) \\ \ga(w) &\de(w)
\end{matrix}\right)}{\det(z)\det^{-1}(z)}
\end{equation*}
and next we use \eqref{def:pairingforhbialgebras} and
\eqref{eq:ellpairingdet}. E.g. for the pairing with $\al(w)$ we
get, using
$\De^{\E^\cop}(\al(w))=\al(w)\otimes\al(w) + \ga(w)\otimes\be(w)$,
\begin{equation*}
\begin{split}
T_{-1} &= \pair{\al(w)}{\det(z)} T_1 \pair{\al(w)}{\det^{-1}(z)} +
\pair{\ga(w)}{\det(z)} T_{-1} \pair{\be(w)}{\det^{-1}(z)} \\
&= q\frac{\theta(w/q^2z)}{\theta(w/z)} \pair{\al(w)}{\det^{-1}(z)},
\end{split}
\end{equation*}
since the second term cancels. This proves the first statement
in
\begin{equation}\label{eq:pairinginvdet}
\begin{split}
&\pair{\left( \begin{matrix} \al(w) &\be(w) \\ \ga(w) &\de(w)
\end{matrix}\right)}{\det^{-1}(z)}=
q^{-1}\frac{\theta(w/z)}{\theta(q^{-2}w/z)}
\left( \begin{matrix} T_{-1} & 0 \\ 0 & T_1
\end{matrix}\right), \\
&\pair{\det^{-1}(w)}{\left( \begin{matrix} \al(z) &\be(z)
\\ \ga(z) &\de(z) \end{matrix}\right)}=
q^{-1}\frac{\theta(q^2w/z)}{\theta(w/z)}
\left( \begin{matrix} T_{-1} & 0 \\ 0 & T_1
\end{matrix}\right),
\end{split}
\end{equation}
the other ones being proved similarly.

\begin{lemma}\label{lem:ellipticpairing}
The pairing on $\E^\cop\times\E\to \Dh$
defined on the generators by \eqref{eq:elldualgenerator}
and \eqref{eq:pairinginvdet} makes $\E^\cop$ and
$\E$ paired as \hHopfster s using
\eqref{eq:pairingXstara}.
\end{lemma}

\begin{proof}
Since the $R$-matrix in \eqref{eq:ellR} is a solution to the
quantum dynamical Yang-Baxter
equation \eqref{eq:QDYBE}, we already obtain the pairing on the
level of \hbialg s by Remark \eqref{rmk:cobraidingandRmatrix}(ii).
A straightforward check shows that
\eqref{eq:pairingforhHopfalgebras} holds on the level of
generators, including $\det^{-1}(z)$, so that we obtain the
pairing on the level of \hHopfalg s. Finally, it remains to be
checked that the $\ast$-structures and the pairing are compatible
using \eqref{eq:pairingXstara}, and using Proposition
\ref{prop:astandpairing} it suffices to check this for the
generators.
\end{proof}

\begin{remark}
For the case of the dynamical quantum group associated
with the rational $R$-matrix \eqref{eq:ratR} Rosengren
\cite[\S 4]{Rose} has calculated the radical of the
pairing, and given it an appropriate representation theoretic
meaning.
For the radical of the pairing for the elliptic $U(2)$
dynamical quantum group this is not known.
\end{remark}

\section{Actions arising from pairings}\label{sec:actionsfrompairings}

In this section we consider two \hbialg s with a pairing.
The pairing can be used to construct a natural action of
one \hbialg\ on the other. This is the content of
Theorem \ref{thm:defaction}, and we also consider the weak
formulation of this action and its intertwining properties
with the \hcoalg -structure.

\begin{thm}\label{thm:defaction}
Let $\U$, $\A$ be \hbialg s with a pairing
$\pair{\cdot}{\cdot} \colon \U\times \A \to \Dh$ as
in Definition \ref{def:pairingforhbialgebras}.
For $X\in \U_{\al\be}$
and $a\in \A$ the elements
\begin{subequations}\label{eq:defaction}
\begin{align}
\label{eq:defadotX}
a\cdot X &= \bigl( T_\al \pair{X}{\cdot}\otimes \id\bigr) \De^\A(a)
= \sum_{(a)} T_\al \pair{X}{a_{(1)}} \otimes a_{(2)}
= \sum_{(a)} \mu_l^\A(T_\al \pair{X}{a_{(1)}}\1) a_{(2)}, \\
\label{eq:defXdota}
X\cdot a &= \bigl( \id \otimes \pair{X}{\cdot}T_\be\bigr) \De^\A(a)
= \sum_{(a)} a_{(1)} \otimes \pair{X}{a_{(2)}}T_\be
= \sum_{(a)} \mu_r^\A(\pair{X}{a_{(2)}}T_\be \1 ) a_{(1)},
\end{align}
\end{subequations}
are well-defined elements in $\A$, and
for $X\in \U_{\al\be}$ and $a\in \A_{\ga\de}$
we have $X\cdot a \in \A_{\ga,\de+\al-\be}$ and
$a\cdot X \in \A_{\ga+\be-\al,\de}$.
Then \eqref{eq:defaction} defines a right and left action of
$\U$ on $\A$;
\begin{equation*}
a\cdot(XY) = (a\cdot X)\cdot Y,\qquad
(XY)\cdot a = X\cdot (Y\cdot a), \qquad
\forall\, a\in \A, \forall\, X,Y \in \U.
\end{equation*}
The left and right action commute;
$(X\cdot a)\cdot Y = X\cdot (a\cdot Y)$, $\forall X,Y\in\U$,
$\forall a\in \A$. Moreover,
with $\De^\U(X)=\sum_{(X)} X_{(1)}\otimes X_{(2)}$,
\begin{equation*}
\begin{split}
X\cdot (ab) &= \sum_{(X)}\bigl( X_{(1)}\cdot a
\bigr)\bigl(X_{(2)}\cdot b\bigr), \quad
(ab)\cdot X = \sum_{(X)}\bigl( a\cdot X_{(1)}
\bigr)\bigl(b\cdot X_{(2)}\bigr),\\
X\cdot 1_\A &= \mu_r^\A(\ep^\U(X)\1 ), \qquad
1_\A\cdot X = \mu^\A_l (T_\al \ep^\U(X)\1 ), \qquad X\in \U_{\al\be}.
\end{split}
\end{equation*}
If $\U$ and $\A$ are paired \hHopfalg s,
then the action of $\U$ on $\A$
satisfies $S^\A(a\cdot S^\U(X))=X\cdot S^\A(a)$ for all
$X\in\U$ and all $a\in\A$.
If $\U$ and $\A$ are paired \hHopfster s with invertible antipodes,
then the action of $\U$ on $\A$ satisfies
$X\cdot a^\ast = \bigl( S^\U(X)^\ast\cdot a\bigr)^\ast
= \bigl( S^\U(X^\da)\cdot a\bigr)^\ast$ and
$a^\ast\cdot X = \bigl((a\cdot S^\U(X)^\ast\bigr)^\ast
= \bigl((a\cdot S^\U(X^\da)\bigr)^\ast$for all
$X\in\U$ and all $a\in\A$.
\end{thm}

The action of Theorem \ref{thm:defaction} can be formulated in terms of the
pairing. Because the pairings are in general not assumed to
be non-degenerate, this is a weaker statement and this is given
in Proposition  \ref{prop:weakaction}. We postpone its proof to the
end of the section.

\begin{prop}\label{prop:weakaction} The action defined in
\eqref{eq:defaction} satisfies
\begin{equation*}
\pair{Y}{X\cdot a} = \pair{YX}{a}\circ T_\be, \qquad
\pair{Y}{a\cdot X} = T_\al \circ \pair{XY}{a}
\end{equation*}
for all $X\in\U_{\al\be}$, $Y\in\U$ and $a\in\A$.
\end{prop}

\begin{proof}[Proof of Theorem \ref{thm:defaction}]
We prove the statements for the right
action of $\U$ on $\A$, the results for the left action being
proved analogously. We put $\De=\De^\A$, $\mu_l=\mu_l^\A$,
$\mu_r=\mu_r^\A$, $S=S^\A$ for the proof.

Since the comultiplication $\De \colon
\A\to  \A\wtt \A$ is well-defined, we only have to show
that $T_\al \pair{X}{\cdot}\otimes \id$ preserves
the relation \eqref{eq:defotimesMh}, i.e. for $f\in\Mh$,
$a,b\in\A$ we need to check
$$
T_\al \pair{X}{\mu_r(f)a}\otimes b =
T_\al \pair{X}{a}\otimes \mu_l(f)b.
$$
Take $X\in\U_{\al\be}$, $a\in \A_{\ga\de}$, then
\eqref{eq:pairing-rightmoment} shows that the
left hand side equals
\begin{equation*}
T_\al \pair{X}{\mu_r(f)a}\otimes b 
= T_\al \pair{X}{a\mu_r(T_\de f)}\otimes b =
T_\al \pair{X}{a}\circ (T_\de f) \otimes b 
=
 \mu_l(T_\al \pair{X}{a} (T_\de f)\1)\, b
\end{equation*}
and the right hand side equals
$\mu_l(T_\al \pair{X}{a}\1) \mu_l(f)b$
and the equality follows from \eqref{eq:trividinDh},
since by \eqref{eq:pairing-target} we
have $T_\al \pair{X}{a}\in (\Dh)_{\de,\be+\ga-\al}$.

Next, for $a\in\A_{\ga\de}$ we write
$\De(a)=\sum_{(a)} a_{(1)}\otimes a_{(2)}$ with
$a_{(1)}\in \A_{\ga\eta}$, $a_{(2)}\in \A_{\eta\de}$.
With $X\in\U_{\al\be}$ we have
$T_\al\pair{X}{a_{(1)}}\in (\Dh)_{\eta,\be+\ga-\al}$,
so this is only non-zero for $\eta=\be+\ga-\al$.
Hence, $a\cdot X\ \in \A_{\be+\ga-\al,\de}$.
So in particular, $\cdot X \colon \A\to \A$
does not preserve the grading, but
$(\mu_r(f)a)\cdot X = \mu_r(f)(a\cdot X)$
is immediate from \eqref{eq:defadotX}
and \eqref{eq:momentAtildetensorB}.
For the left moment map and $X\in\U_{\al\be}$ we have
\begin{equation}\label{eq:leftmomentrightaction}
\begin{split}
(\mu_l(f)a)\cdot X &=
\sum_{(a)}
\mu_l(T_\al\pair{X}{\mu_l(f)a_{(1)}}\1) a_{(2)} =
\sum_{(a)}
\mu_l(T_\al \circ f\circ \pair{X}{a_{(1)}}\1) a_{(2)}
\\ &=
\mu_l(T_\al f) \sum_{(a)}
\mu_l(T_\al \pair{X}{a_{(1)}}\1) a_{(2)} =
\mu_l(T_\al f) (a\cdot X)
\end{split}
\end{equation}
using
\eqref{eq:momentAtildetensorB},  \eqref{eq:defadotX} and
\eqref{eq:pairing-leftmoment}.

To show that this defines an action we write $\De(a)$ for $a\in\A_{\ga\de}$
as above and we use the coassociativity
\begin{equation*}
\begin{split}
(\id\otimes\De)\De(a) = (\De\otimes\id)\De(a) &=
\sum_{(a)} a_{(1)}\otimes a_{(2)}\otimes a_{(3)},
\quad a_{(1)}\in \A_{\ga\eta}, \ a_{(2)}\in \A_{\eta\rho}, \
a_{(3)}\in \A_{\rho\de}.
\end{split}
\end{equation*}
With this convention we have for $X\in \U_{\al\be}$, $Y\in
\U_{\si\tau}$, and hence $XY\in \U_{\al+\si,\be+\tau}$,
\begin{equation*}
\begin{split}
a\cdot (XY) &= \sum_{(a)} T_{\al+\si}\pair{XY}{a_{(1)}}
\otimes a_{(2)} =
\sum_{(a)} T_{\al+\si}\pair{X}{a_{(1)}} T_\eta \pair{Y}{a_{(2)}}
\otimes a_{(3)} \\ &=
\sum_{(a)} \mu_l(T_{\al+\si}\pair{X}{a_{(1)}}
T_\eta \pair{Y}{a_{(2)}} \1) a_{(3)}
\end{split}
\end{equation*}
using \eqref{eq:pairing-leftproduct}.
On the other hand, using the $\C$-linearity of the action
and \eqref{eq:leftmomentrightaction} we have
\begin{equation*}
\begin{split}
(a\cdot X)\cdot Y &=
\Bigl(  \sum_{(a)} \mu_l(T_\al \pair{X}{a_{(1)}}\1) a_{(2)}
\Bigr)\cdot Y =
\sum_{(a)} \mu_l(T_\si
T_\al \pair{X}{a_{(1)}}\1)
\bigl(a_{(2)}\cdot Y\bigr) \\
 &= \sum_{(a)} \mu_l(T_{\si+\al} \pair{X}{a_{(1)}}\1)
\mu_l(T_\si \pair{Y}{a_{(2)}}\1) a_{(3)}
\end{split}
\end{equation*}
and since
$T_{\si+\al} \pair{X}{a_{(1)}} \in (\Dh)_{\eta-\si,\ga+\be-\si-\al}$,
we see that $a\cdot (XY)=(a\cdot X)\cdot Y$
follows from \eqref{eq:trividinDh} and $(\Dh)_{\al\be}=\{0\}$
for $\al\not=\be$.

For the commutativity of the left and right action we
take $X\in\U_{\al\be}$,
$Y\in\U_{\si\tau}$, and using \eqref{eq:actionmomentA},
\begin{equation*}
\begin{split}
(X\cdot a)\cdot Y &=
\Bigl( \sum_{(a)} \mu_r(\pair{X}{a_{(2)}}T_\be \1
)a_{(1)}\Bigr)\cdot Y =
\sum_{(a)} \mu_r(\pair{X}{a_{(2)}}T_\be \1)
\bigl( a_{(1)} \cdot Y\bigr) \\
&= \sum_{(a)} \mu_r(\pair{X}{a_{(3)}}T_\be \1)
\mu_l(T_\si\pair{Y}{a_{(1)}}\1) a_{(2)}
\end{split}
\end{equation*}
and similarly
\begin{equation*}
\begin{split}
X\cdot (a\cdot Y) &= X\cdot \Bigl( \sum_{(a)}
\mu_l(T_\si\pair{Y}{a_{(1)}}\1) a_{(2)}\Bigr) = \sum_{(a)}
\mu_l(T_\si\pair{Y}{a_{(1)}}\1) \bigl( X\cdot a_{(2)})
\\ &=
\sum_{(a)}
\mu_l(T_\si\pair{Y}{a_{(1)}}\1)
\mu_r(\pair{X}{a_{(3)}}T_\be \1)a_{(2)},
\end{split}
\end{equation*}
which proves the statement.

For $X\in\U_{\al\be}$ we see, using the \halg\ homomorphism property
of $\De$ and \eqref{eq:pairing-rightproduct},
\begin{equation*}
\begin{split}
(ab)\cdot X &= \sum_{(a)}\sum_{(b)} \mu_l(T_\al
\pair{X}{a_{(1)}b_{(1)}}\1 )
a_{(2)}b_{(2)}\\
&= \sum_{(a)}\sum_{(b)} \mu_l(T_\al
\sum_{(X)}\pair{X_{(1)}}{a_{(1)}}T_\eta\pair{X_{(2)}}{b_{(1)}}\1 )
a_{(2)}b_{(2)}
\end{split}
\end{equation*}
using $\De(X)=\sum_{(X)} X_{(1)}\otimes X_{(2)}$ with
$X_{(1)}\in\U_{\al\eta}$, $X_{(2)}\in\U_{\eta\be}$. Writing $\De
(a) = \sum_{(a)} a_{(1)}\otimes a_{(2)}$, $a_{(1)}\in
\A_{\ga\rho}$, $a_{(2)}\in \A_{\rho\de}$ we see that
$T_\al\pair{X_{(1)}}{a_{(1)}} \in(\Dh)_{\rho,\ga+\eta-\al}$ and by
\eqref{eq:trividinDh},
$$
\mu_l(T_\al \pair{X_{(1)}}{a_{(1)}}T_\eta\pair{X_{(2)}}{b_{(1)}}\1
) = \mu_l(T_\al \pair{X_{(1)}}{a_{(1)}}\1)
\mu_l(T_{\al-\ga}\pair{X_{(2)}}{b_{(1)}}\1 ).
$$
This gives
\begin{equation*}
\begin{split}
(ab)\cdot X &= \sum_{(a)}\sum_{(b)}\sum_{(X)} \mu_l(T_\al
\pair{X_{(1)}}{a_{(1)}}\1)
\mu_l(T_{\al-\ga}\pair{X_{(2)}}{b_{(1)}}\1 )
a_{(2)}b_{(2)} \\
&= \sum_{(a)}\sum_{(b)}\sum_{(X)} \mu_l(T_\al
\pair{X_{(1)}}{a_{(1)}}\1) \, a_{(2)}
\mu_l(T_{\rho+\al-\ga}\pair{X_{(2)}}{b_{(1)}}\1 ) b_{(2)}
\end{split}
\end{equation*}
and since $(\Dh)_{\rho,\ga+\eta-\al}=\{0\}$ unless
$\rho+\al-\ga=\eta$ we find
\begin{equation*}
\begin{split}
(ab)\cdot X &=\sum_{(X)}
\Bigl( \sum_{(a)}\mu_l(T_\al
\pair{X_{(1)}}{a_{(1)}}\1)a_{(2)}\Bigr)
\Bigl(\sum_{(b)}\mu_l(T_\eta\pair{X_{(2)}}{b_{(1)}}\1 )
b_{(2)}\Bigr)\\
&= \sum_{(X)} \bigl( a\cdot X_{(1)}\bigr) \bigl( b\cdot X_{(2)}\bigr).
\end{split}
\end{equation*}
This proves the statements for the right action.

Assume now that $\U$ and $\A$ are paired \hHopfalg s. For
the proof of $X\cdot S(a) = S(a\cdot S^\U(X))$ we use
$\A\otimes\Dh\cong\A\cong \Dh\otimes \A$, see
\S \ref{ssec:algebraicnotions}.
We take
$X\in\U_{\al\be}$, so
\begin{equation*}
\begin{split}
X\cdot S(a) &= (\id\otimes \pair{X}{\cdot}T_\be)\De(S(a))
=  (\id\otimes \pair{X}{\cdot}T_\be) \circ P \circ
(S\otimes S)\De(a) \\
&= P\circ ( \pair{X}{\cdot}T_\be \otimes \id) \circ
(S\otimes S)\De(a)
=\sum_{(a)} S(a_{(2)}) \otimes \pair{X}{S(a_{(1)})}T_\be \\
&= \sum_{(a)} S(a_{(2)}) \otimes S^\Dh(T_{-\be} \pair{S^\U(X)}{a_{(1)}})
= S\bigl( T_{-\be} \pair{S^\U(X)}{a_{(1)}}\otimes a_{(2)}\bigr)
= S( a\cdot S^\U(X)),
\end{split}
\end{equation*}
since $S^\U(X)\in \U_{-\be,-\al}$.

Finally, if $\U$ and $\A$ are in paired \hHopfster s,
we have in the same way for $X\in\U_{\al\be}$
\begin{equation*}
\begin{split}
X\cdot a^\ast &= \sum_{(a)} a_{(1)}^\ast \otimes
\pair{X}{a_{(2)}^\ast}T_\be =
\sum_{(a)} a_{(1)}^\ast \otimes
T_{-\al} \bigl(\pair{S^\U(X)^\ast}{a_{(2)}}\bigr)^\ast \\
&=  \Bigl( \sum_{(a)} a_{(1)} \otimes
\pair{S^\U(X)^\ast}{a_{(2)}} T_{\bar\al} \Bigr)^\ast
= \bigl( S^\U(X)^\ast \cdot a\bigr)^\ast,
\end{split}
\end{equation*}
since $S^\U(X)^\ast \in \U_{\bar\be\bar\al}$.
\end{proof}

For later use we note, see also the proof of Theorem
\ref{thm:defaction},  that
for $X\in\U_{\al\be}$, $a\in\A_{\ga\de}$ we have
\begin{equation}\label{eq:actionmomentA}
\begin{array}{ll}
X\cdot (\mu_l^\A(f)a) = \mu_l^\A(f)(X\cdot a),  &
X\cdot (\mu_r^\A(f)a) = \mu_r^\A(T_{-\al}f) (X\cdot a), \\ [1mm]
(\mu_l^\A(f)a)\cdot X = \mu_l^\A(T_\al f) (a\cdot X), &
(\mu_r^\A(f)a)\cdot X = \mu_r^\A(f) (a\cdot X),
\end{array}
\end{equation}
and
\begin{equation}\label{eq:actionmomentU}
\begin{array}{ll}
(\mu_l^\U(f)X)\cdot a = \mu_r^\A(f)(X\cdot a),  &
(\mu_r^\U(f)X)\cdot a = \mu_r^\A(T_{\be-\al-\de}f) (X\cdot a), \\ [1mm]
a\cdot (\mu_l^\U(f)X) = \mu_l^\A(T_\al f) (a\cdot X), &
a\cdot (\mu_r^\U(f)X) = \mu_l^\A(T_{\al-\ga}f) (a\cdot X).
\end{array}
\end{equation}
In particular, taking $X=1\in \U_{00}$ in
\eqref{eq:actionmomentU} gives the left
and right action of the left and right moment map on $\A$.

The following Proposition describes the interaction of the
action and the \hcoalg\ structure of $\A$, and this is
useful in constructing corepresentations using
invariance properties in terms of the action defined
in Theorem \ref{thm:defaction}.

\begin{prop}\label{prop:actionDeltaepsilon}
With the notation and assumptions as in Theorem \ref{thm:defaction} we have
\begin{equation*}
\begin{matrix}
\De^\A(X\cdot a) = (\id\otimes X\cdot) \De^\A(a), \quad &
\De^\A(a\cdot X) = (\cdot X \otimes\id)  \De^\A(a), \\[2mm]
\ep^\A(X\cdot a) = \pair{X}{a}\circ T_\be, &
\ep^A(a\cdot X) = T_\al\circ \pair{X}{a},
\end{matrix}
\end{equation*}
for $X\in\U_{\al\be}$ and $a\in\A$.
\end{prop}

\begin{proof} By \eqref{eq:actionmomentA} it follows
that $a\otimes (X\cdot\mu_l^\A(f)b) =
a\otimes \mu_l^\A(f) (X\cdot b) = \mu_r^\A(f)a\otimes
(X\cdot b)$ so that $\id\otimes X\cdot \colon \A\wtt\A \to
\A\wtt\A$ is well-defined, cf. \eqref{eq:defotimesMh}.
Now for $X\in\U_{\al\be}$
\begin{equation*}
\begin{split}
\De^\A(X\cdot a) &= \sum_{(a)}
\De^\A\Bigl( \mu_r^\A( \pair{X}{a_{(2)}}T_\be\1 )a_{(1)}\Bigr)
= \sum_{(a)} \bigl( 1\otimes \mu_r^\A( \pair{X}{a_{(2)}}T_\be\1 )\bigr)
\De^\A(a_{(1)}) \\
&= \sum_{(a)} a_{(1)}\otimes \mu_r^\A( \pair{X}{a_{(3)}}T_\be\1 )a_{(2)}
= \sum_{(a)} a_{(1)}\otimes (X\cdot a_{(2)})
\end{split}
\end{equation*}
using $\De^\A(a)=\sum_{(a)}a_{(1)} \otimes a_{(2)}$,
$(\De^\A\otimes\id)\De^\A(a)=(\id\otimes\De^\A)\De^\A(a)=
\sum_{(a)}a_{(1)} \otimes a_{(2)}\otimes a_{(3)}$, so  the
result is a consequence of the coassociativity.

Now for $X\in\U_{\al\be}$
\begin{equation*}
\ep^\A(X\cdot a) = (\ep^\A \otimes \pair{X}{\cdot}T_{\be})\De^\A(a)
= (\id_\Dh \otimes \pair{X}{\cdot}T_{\be})(\ep^\A \otimes \id_\A)\De^\A(a)
= \pair{X}{a}\circ T_\be
\end{equation*}
using the counit axiom and the appropriate identifications as
discussed in \S \ref{ssec:algebraicnotions}.
The statements for the right action are proved similarly.
\end{proof}

It remains to prove Proposition \ref{prop:weakaction}.

\begin{proof}[Proof of Proposition \ref{prop:weakaction}]
We prove the first statement, the other one being proved
analogously. Assume $X\in\U_{\al\be}$,
$Y\in\U_{\si\tau}$,
$a\in\A_{\ga\de}$,
$\De^\A(a) = \sum_{(a)} a_{(1)}\otimes a_{(2)}$
with $a_{(1)}\in\A_{\ga\rho}$, $a_{(2)}\in\A_{\rho\de}$,
\begin{equation*}
\begin{split}
\pair{Y}{X\cdot a} &= \pair{Y}{\sum_{(a)}
\mu_r^\A(\pair{X}{a_{(2)}}T_\be\1)a_{(1)}}
= \pair{Y}{\sum_{(a)} a_{(1)} \mu_r^\A(T_\rho \pair{X}{a_{(2)}}T_\be\1)} \\
&=  \sum_{(a)} \pair{Y}{ a_{(1)}} \circ\bigl( T_\rho
\pair{X}{a_{(2)}}T_\be\1\bigr).
\end{split}
\end{equation*}
Write $\pair{X}{a_{(2)}} =fT_{-\be-\rho} \in
\bigl( \Dh\bigr)_{\al+\de,\be+\rho}$,
$\pair{Y}{a_{(1)}} =gT_{-\si-\rho} \in
\bigl( \Dh\bigr)_{\si+\rho,\tau+\ga}$, so that the summand
can be rewritten as
$g\circ (T_{-\si}f) \circ T_{-\si-\rho} \in \Dh$.
On the other hand, it is straightforward that this is also
equal to
$\pair{Y}{a_{(1)}}\circ T_\rho\circ \pair{X}{a_{(2)}} \circ T_\be$,
so that
\begin{equation*}
\pair{Y}{X\cdot a} = \sum_{(a)} \pair{Y}{ a_{(1)}} \circ T_\rho
\circ \pair{X}{a_{(2)}} \circ T_\be = \pair{YX}{a}\circ T_\be.
\end{equation*}
\end{proof}

\begin{remark}\label{rmk:nodouble}
Having paired \hHopfalg s $\U$ and $\A$ one might expect to
have an analogue of the double construction of Drinfeld, see
Lu \cite{Lu}  for the case of Hopf algebroids. However,
if we want to mimic the classical construction we see that
for elements $X,Y\in\U$, with $\De^\U(X)=\sum_{(X)}X_{(1)}\otimes X_{(2)}$,
and $a\in \A$ neither the element
$\sum_{(X)} S^\U(X_{(2)})\cdot a \cdot X_{(1)}$ nor
$\sum_{(X)} X_{(1)}Y S^\U(X_{(2)})$, or any other similarly defined element,
is well-defined, i.e. respects the relation \eqref{eq:defotimesMh} in $\U\wtt\U$.
For a double construction for $\times_R$-bialgebras
using a different pairing, see \cite[\S 6]{Scha}.
\end{remark}

\section{Representation and corepresentation theory}\label{sec:repres}

Given two \hbialg s with a pairing, we consider
the relation between dynamical representations of the one
\hbialg\ and the corepresentations of the other \hbialg .
For \hHopfster s we rephrase the notion of unitarisability
introduced in \cite{KoelR} and consider it in the context of two
\hHopfster s equipped with a pairing, and we show that this
extends the previously introduced notion of unitarisability. For the
elliptic $U(2)$ dynamical quantum group we make this
explicit, and this enables us to calculate the pairing
between matrix elements of irreducible corepresentations
in terms of elliptic hypergeometric series.
We start by recalling the notions of corepresentations and
dynamical representations.

A {\em \hspa}\ $V$ is a vector space over $\Mh$, which is a also
a diagonalisable $\h$-module, $V=\bigoplus_{\al\in\h^\ast} V_\al$,
$V_\al = \{ v\in V\mid H\cdot v=\al(H)v,\, \forall H\in\h \}$,
with $\Mh V_\al\subseteq V_\al$ for all $\al$.
A {\em morphism} of \hspa s is a $\Mh$-linear $\h$-invariant map.
In case we want to emphasize the multiplication of $v\in V$
with $f\in\Mh$ we also write $\mu_V(f)v= fv$.
Note that the definition of a \hspa\ implies that $V$ has a
weight space decomposition, and we can speak of a homogeneous
vector $v\in V_\al$ as a vector of weight $\al\in\hs$.

For a \halg\ $\A$ and a \hspa\ $V$ we define
$\A \wtt V =\bigoplus_{\al\be} \A_{\al\be}\otimes_\Mh V_\be$,
where we mod out by the relations
$\mu_r^\A(f)a\otimes v =a\otimes fv$ for all $f\in\Mh$,
$a\in \A$, $v\in V$, cf. \eqref{eq:defotimesMh}.
We make $\A \wtt V$ into a \hspa\ by the grading
$\A_{\al\be}\otimes_\Mh V_\be \subseteq (\A\wtt V)_\al$
for any $\be\in\hs$, and by $\mu_{\A \wtt V}(f)(a\otimes v) =\mu_l^\A(f)a \otimes v$.
This is easily checked. This definition is compatible with the
tensor product of \halg s, see \eqref{eq:defAtildetensorB},
in the sense $(\A\wtt\B)\wtt V = \A\wtt(\B\wtt V)$
as \hspa s for \halg s $\A$ and $\B$ and a \hspa\ $V$.

Similarly we define
$V \wtt \A =\bigoplus_{\al\be} V_\al\otimes_\Mh \A_{\al\be}$,
where we mod out by the relations
$v\otimes \mu_l^\A(f)a =fv\otimes a$ for all $f\in\Mh$,
$a\in \A$, $v\in V$, cf. \eqref{eq:defotimesMh}.
Then  $V \wtt \A$ is a \hspa\ with grading
$V_\al\otimes_\Mh \A_{\al\be} \subseteq (\A\wtt V)_\be$
for all $\al$, and $\mu_{V\wtt\A}(f)(v\otimes a) =v\otimes \mu_r^\A(f)a$.
This time we have the compatibility
relation $V\wtt (\A\wtt\B) = (V\wtt\A)\wtt\B$
as \hspa s for \halg s $\A$ and $\B$ and a \hspa\ $V$.

In the special case $\A=\Dh$ we obtain
$V\wtt\Dh\cong V \cong \Dh\wtt V$ as \hspa s with the
isomorphism given by
$v\otimes fT_{-\al}\cong fv\cong fT_{-\al}\otimes v$
for $f\in\Mh$, $v\in V_\al$.

Now a {\em left corepresentation} of a \hcoalg\ $\A$ on a \hspa\
$V$ is a \hspa \ morphism $L\colon V\to \A\wtt V$ satisfying
$(\De^\A\otimes\id)\circ L = (\id\otimes L)\circ L$ and
$(\ep^\A\otimes\id)\circ L =\id$. Note that for the first
requirement we use $(\A\wtt\A)\wtt V = \A\wtt(\A\wtt V)$ and for
the second we use $V \cong \Dh\wtt V$. Similarly, a {\em right
corepresentation}  of a \hcoalg\ $\A$ on a \hspa\ $V$ is a \hspa
-morphism $R\colon V\to V\wtt\A$ satisfying $(R\otimes\id)\circ R
= (\id\otimes \De^\A)\circ R$ and $(\id\otimes\ep^\A)\circ R
=\id$.

An {\em intertwiner} of two left, respectively right,
corepresentations $L_1$ and $L_2$, respectively $R_1$ and $R_2$,
in $V_1$ and $V_2$ is a
\hspa\ morphism $\Phi\colon V_1\to V_2$ such that
$L_2\circ \Phi = (\id\otimes \Phi) \circ L_1$,
respectively $R_2\circ \Phi = (\Phi\otimes \id) \circ R_1$.
Note that $\id\otimes \Phi$, respectively $\Phi\otimes \id$,
do factor to a map on $\A\wtt V_1$, respectively $V_1\wtt\A$.

To a \hspa\ $V$ we can associate a \halg\ $D_{\hs,V}$
as follows. First, we define the graded
subspace $\bigl(D_{\hs,V}\bigr)_{\al\be}$
as the space of $\C$-linear operators $U\colon V\to V$
satisfying $U(fv)=(T_{-\be}f)(Uv)$ for all $f\in\Mh$, $v\in V$,
and $U(V_\ga)\subseteq V_{\ga-\al+\be}$.
The moment maps
$\mu_l,\mu_r\colon \Mh\to \bigl(D_{\h,V}\bigr)_{00}$
are defined by
$\mu_r(f)v= fv$, $\mu_l(f)v= (T_{-\ga}f)v$ for $v\in V_\ga$.
One easily checks that this makes
$D_{\hs,V}=\bigoplus_{\al,\be\in\hs}
\bigl(D_{\hs,V}\bigr)_{\al\be}$
into a \halg .
A {\em dynamical representation} of a \halg\ $\A$ on a \hspa\ $V$
is a \halg-morphism $\pi\colon \A\to D_{\hs,V}$.
For $\pi_i\colon \A\to D_{\hs,V_i}$, $i=1,2$,
dynamical representations
of a \halg\ $\A$ we say that a \hspa\ morphim
$\Phi\colon V_1\to V_2$ is
an {\em intertwiner} of $\pi_1$ and $\pi_2$ if
$\Phi\circ \pi_1(a) = \pi_2(a)\circ\Phi$ for all $a\in\A$.

Recall from Remark \ref{rmk:pairingforhbialgebras}
the pairing between
\halg s and \hcoalg s. The applications of
Proposition \ref{prop:cotodynrep} are for \hbialg s,
we formulate it in the slightly more general
case.

\begin{prop}\label{prop:cotodynrep}
Let $\U$ be a \halg\ and $\A$ be \hcoalg\  equipped with
a pairing, and
let $V$ be a \hspa .
\par\noindent {\rm (i)}
Let $R\colon V\to V\wtt \A$ be a right corepresentation
of the \hcoalg\ $\A$, then
$\pi(X)v = (\id\otimes \pair{X}{\cdot}T_\be)R(v)$ for
$X\in\U_{\al\be}$,
defines a \halg\ morphism $\pi\colon \U \to (D_{\hs,V})^\lr$
of $\U$ on $V$, hence $\pi\colon \U^\lr \to D_{\hs,V}$
defines a dynamical representation of $\U^\lr$ on $V$.
\par\noindent {\rm (ii)}
Let $L\colon V\to \A\wtt V$ be a left corepresentation
of the \hcoalg\ $\A$, then
$\pi(X)v = (T_{\al}\pair{X}{\cdot}\otimes\id)L(v)$ for
$X\in\U_{\al\be}$, defines a \halg\ homomorphism
$\pi\colon \U^\opp \to (D_{\hs,V})^\lr$. In particular,
$\pi\colon (\U^\opp)^\lr \to D_{\hs,V}$ defines a
dynamical representation of $(\U^\opp)^\lr$ on $V$.
Moreover, if $\U$ is \hHopfalg , then
$X\mapsto \pi(S^\U(X))$ defines a dynamical representation
of $\U$ on $V$.
\end{prop}

In case $\A$ is equipped with a cobraiding,
see Definition \ref{defn:cobraiding}, we can use
$\A^\lr=\A^\cop$ as \halg s. Also
$(\A^\cop)^\opp = (\A^\opp)^\cop$, and
this immediately implies the next result.

\begin{cor}\label{cor:propcotodynrep}
Let $\A$ be a \hbialg\ equipped with a cobraiding.
\par\noindent {\rm (i)} A right corepresentation of $\A$
on $V$
gives rise to a dynamical representation of $\A$ on $V$.
\par\noindent {\rm (ii)} A left corepresentation of $\A$
on $V$
gives rise to a dynamical representation of $\A^\opp$ on $V$.
\end{cor}

\begin{proof}[Proof of Proposition \ref{prop:cotodynrep}]
The proof is
a minor variation on the proof of Theorem
\ref{thm:defaction}, so we do not give all details.

For $v\in V_\ga$ write $Rv=\sum_{(v)} v_{(1)}\otimes a_{(2)}$
with $v_{(1)}\in V_\de$, $a_{(2)}\in\A_{\de\ga}$, so that
for $X\in\U_{\al\be}$
\begin{equation*}
\pi(X) v= (\id\otimes \pair{X}{\cdot}T_\be) R(v)
= \sum_{(v)} \mu_V(\pair{X}{a_{(2)}}T_\be\1)\, v_{(1)},
\end{equation*}
and since $\pair{X}{a_{(2)}}T_\be\in (\Dh)_{\al+\ga-\be,\de}$ we
see that $\pi(X)v$ is well-defined and that we only find a
contribution to the sum for $\de=\ga+\al-\be$. So $\pi(X)
V_\ga\subseteq V_{\ga+\al-\be}$. Next, for $f\in\Mh$ we find
$R(fv)= f (Rv) = \sum_{(v)} v_{(1)}\otimes \mu_r^\A(f)a_{(2)}$, so
that
\begin{equation*}
\begin{split}
\pi(X) (fv) &= \sum_{(v)}
\mu_V(\pair{X}{\mu_r^\A(f)a_{(2)}}T_\be\1)\, v_{(1)} =
\sum_{(v)}
\mu_V(\pair{X}{a_{(2)}}\, (T_\ga f)T_\be\1)\, v_{(1)}
\\ &= (T_{-\al} f)
\sum_{(v)}
\mu_V(\pair{X}{a_{(2)}}\, T_\be\1)\, v_{(1)}
= (T_{-\al} f)\,  \pi(X)v
\end{split}
\end{equation*}
using \eqref{eq:pairing-rightmoment} and
$\pair{X}{a_{(2)}}\in (\Dh)_{\al+\ga,\de+\be}$. This also
shows the $\C$-linearity, and hence
$\pi(X)\in (D_{\hs,V})_{\be\al}$. Furthermore, $\pi\colon\U
\to D_{\hs,V}$ is $\C$-linear, and the relation
$\pi(XY)=\pi(X)\pi(Y)$ is proved as the similar statement
for the left action in Theorem \ref{thm:defaction}.
For the actions of the left and right moment maps we
find
\begin{equation*}
\begin{split}
\pi(\mu^\U_l(f))v &= \sum_{(v)}
\mu_V(\pair{\mu_l^\U(f)}{a_{(2)}}T_\be\1)\, v_{(1)} =
\mu_V(f \ep^\A(a_{(2)})T_\be\1)\, v_{(1)}
\\ &= \mu_V(f) (\id\otimes\ep^\A)\, Rv = \mu_V(f)\, v, \\
\pi(\mu^\U_r(f))v &= \sum_{(v)}
\mu_V(\pair{\mu_r^\U(f)}{a_{(2)}}T_\be\1)\, v_{(1)} =
\mu_V( \ep^\A(a_{(2)})\, f\, T_\be\1)\, v_{(1)} \\
&= \mu_V(T_{-\ga}f) (\id\otimes\ep^\A)\, Rv =
\mu_V(T_{-\ga}f)\, v,
\end{split}
\end{equation*}
using \eqref{eq:pairing-leftmoment}, \eqref{eq:pairing-counit-unit}
and $\ep^\A(a_{(2)})\in (\Dh)_{\de\ga}$.
Hence, $\pi\colon \U\to (D_{\hs,V})^\lr$ is
a \halg-homomorpism, and consequently
$\pi\colon \U^\lr\to D_{\hs,V}$ defines a
dynamical representation of $\U^\lr$ on $V$.

For the second statement we put for $v\in V_\ga$
$Lv=\sum_{(v)} a_{(1)} \otimes v_{(2)}$ with
$a_{(1)}\in \A_{\ga,\de}$, $v_{(2)}\in V_\de$, so that
for $X\in \U_{\al\be}$ we have
$\pi(X)v= \sum_{(v)} \mu_V(T_\al \pair{X}{a_{(1)}}) v_{(2)}$.
This implies $\pi(X)V_\ga\subseteq V_{\ga+\be-\al}$ and
$\pi(X)[fv] = \mu_V(T_\al f)\pi(X)v$, and hence
$\pi(X)\in (D_{\hs,V})_{-\be,-\al}$ for $X\in\U_{\al\be}$.
We also obtain
$\pi(\mu_l^\U(f))\, v= \mu_V(f)\, v$,
$\pi(\mu_r^\U(f))\, v= \mu_V(T_{-\ga}f)\, v$
using \eqref{eq:pairing-leftmoment},
\eqref{eq:pairing-rightmoment},
\eqref{eq:pairing-counit-unit} and $(\ep\otimes\id)\circ L=\id$.
The proof of $\pi(XY)=\pi(Y)\pi(X)$ proceeds as the
corrresponding statement for the right action in Theorem
\ref{thm:defaction}.
Since $\pi((\U^\opp)_{\al\be})=\pi(\U_{-\al,-\be})
\subseteq (D_{\hs,V})_{\be\al}=(D_{\hs,V})^\lr_{\al\be}$
we obtain the first statement in (ii), and the second
follows immediately.

In case $\U$ is \hHopfalg , put $\tilde\pi(X)=\pi(S^\U(X))$.
Then the antimultiplicativity
of $S$ makes $\tilde\pi$ an algebra (over $\C$) homomorphism
satisfying $\tilde\pi(\U_{\al\be}) \subseteq
(D_{\hs,V})_{\al\be}$, $\tilde\pi(\mu_l^\U(f))v=
\mu_V(T_{-\ga}f)\, v$, $\tilde\pi(\mu_r^\U(f))v= \mu_V(f)v$,
since $S$ interchanges the left and right moment map.
\end{proof}

The dynamical representations of $\U$  and corepresentations
of $\A$ for paired $\U$ and $\A$  are related as
described in the following lemma.

\begin{lemma}\label{lemma:interactioncorep-rep}
Assume the conditions of Proposition \ref{prop:cotodynrep}.
\par\noindent {\rm (i)}
$R\circ \pi(X) = (\id\otimes X\cdot\, ) \circ R$ with $R$ and $\pi$
as in  Proposition \ref{prop:cotodynrep}(i).
\par\noindent {\rm (ii)}
$L\circ \pi(X) = (\, \cdot X\otimes \id) \circ L$ with $L$ and $\pi$
as in  Proposition \ref{prop:cotodynrep}(ii).
\end{lemma}

It follows from \eqref{eq:actionmomentA}
that $\id\otimes X\cdot\,$,
respectively $\, \cdot X\otimes \id$, factors to a map
on $V\wtt\A$, respectively $\A\wtt V$. The
proof of Lemma \ref{lemma:interactioncorep-rep}
is analogous to the proof of Proposition
\ref{prop:actionDeltaepsilon}, and we leave the proof to
the reader.


Next we want to consider unitary corepresentations of a \hbialg\
$\A$ equipped with a $\ast$-operator. This means that $\hs$ and
$\Mh$ are equipped with a complex conjugation. For a \hspa\ $V$ a
{\em form} is a $\Mh$-sesquilinear form $\langle\cdot,\cdot\rangle
\colon V\times V\to \Mh$ satisfying $\pair{v}{v}\not= 0$ in $\Mh$
and $\pair{v}{w}=\overline{\pair{w}{v}}$ for all $v,w\in V$. We
extend this form to $\A\wtt V$ by
\begin{equation}\label{eq:extpairingtoAotV}
\langle\cdot,\cdot\rangle \colon \A\wtt V
\times \A\wtt V \to \A,
\qquad
\pair{a\otimes v}{b\otimes w}= b^\ast \mu_r^\A(\pair{v}{w}) a,
\end{equation}
and to $V\wtt \A$ by
\begin{equation}\label{eq:extpairingtoVotA}
\langle\cdot,\cdot\rangle \colon V\wtt \A
\times V\wtt\A \to \A,
\qquad
\pair{v\otimes a}{w\otimes b}= b^\ast \mu_l^\A(\pair{v}{w}) a.
\end{equation}
We have to check that \eqref{eq:extpairingtoAotV}
and \eqref{eq:extpairingtoVotA} are well-defined.
For \eqref{eq:extpairingtoVotA} we have
$\pair{fv\otimes a}{w\otimes b} =
b^\ast \mu_l^\A(\pair{fv}{w}) a =
b^\ast \mu_l^\A(f) \mu_l^\A(\pair{v}{w}) a$
and $\pair{v\otimes \mu_l^\A(f)a}{w\otimes b}=
b^\ast \mu_l^\A(\pair{v}{w}) \mu_l^\A(f)a$, so
it is well-defined in the first argument.
Similarly
$\pair{v\otimes a}{fw\otimes b} =
b^\ast \mu_l^\A(\pair{v}{fw}) a$ is equal to
$\pair{v\otimes a}{w\otimes \mu_l^\A(f)b} =
b^\ast \mu_l^\A (\bar f)\mu_l^\A(\pair{v}{w}) a$
since the form is sesquilinear. For
the well-definedness of
\eqref{eq:extpairingtoAotV} an analogous argument
can be given.

The extended forms on $\A\wtt V$ and $V\wtt\A$ are no longer
sesquilinear over $\Mh$. It follows directly that
$\pair{v\otimes a}{w\otimes b} = \pair{w\otimes b}{v\otimes
a}^\ast$,  and
$\pair{a\otimes v}{b\otimes w} = \pair{b\otimes w}{a\otimes
v}^\ast$.
For the forms on $\A\wtt V$ we find for $a\otimes v\in (\A\wtt
V)_\ga$, respectively $v\otimes a \in (V\wtt\A)_\ga$,
\begin{equation}\label{eq:extformnonsesquilinear}
\begin{split}
\pair{f(a\otimes v)}{b\otimes w} &=
\pair{a\otimes v}{b\otimes w} \mu_l^\A (T_\ga f), \\
\pair{f(v\otimes a)}{w\otimes b} &= \pair{v\otimes a}{w\otimes b}
\mu_r^\A (T_\ga f).
\end{split}
\end{equation}

Note that in particular such a form
is extended to $\langle\cdot,\cdot \rangle_D
\colon V\times V\to\Dh$
by applying the above procedure to $\A=\Dh$ and using
the identifications $V\wtt\Dh\cong V\cong \Dh\wtt V$.
Since the left and right moment map for $\Dh$ coincide
we find $\langle\cdot,\cdot \rangle_D \colon V\times V\to\Dh$
is given by $\pair{v}{w}_D = T_{\bar\tau}\circ
\pair{v}{w} \circ T_{-\si}$
for $v\in V_\si$, $w\in V_\tau$, where the original
$\pair{v}{w}\in\Mh$ is considered as multiplication
operator. Using this we can rewrite
\eqref{eq:extpairingtoAotV} and
\eqref{eq:extpairingtoVotA} as
follows
\begin{equation}\label{eq:diffformextendedinprodVotA}
\pair{a\otimes v}{b\otimes w}= b^\ast a \otimes \pair{v}{w}_D,
\qquad
\pair{v\otimes a}{w\otimes b}=\pair{v}{w}_D \otimes b^\ast  a,
\end{equation}
using the identification \eqref{eq:identDhotimesAwithA}.

\begin{defn} \label{def:unitarisablecorep}
Let $\A$ be a \hHopfster\ and V a \hspa .
A right corepresentation $R\colon V\to V\wtt\A$ is
unitarisable if there exists a form
$\langle\cdot,\cdot \rangle \colon V\times V
\to \Mh$ such that
$\pair{Rv}{Rw} = \mu_r^\A( \pair{v}{w}_D \1 )$
for all $v,w\in V$.
A left corepresentation  $L\colon V\to \A\wtt V$
is unitarisable if there exists a form
$\langle\cdot,\cdot \rangle \colon V\times V
\to \Mh$ such that
$\pair{Lv}{Lw} =  \mu_l^\A( \pair{v}{w}_D\1 )$
for all $v,w\in V$.
\end{defn}

Unitarisable corepresentations have been introduced
before in \cite[Def.~3.11]{KoelR}, and to see
its relation to
Definition \ref{def:unitarisablecorep},
we choose a homogeneous basis $\{ v_i\}_{i\in I}$
of $V$ as a vector space over $\Mh$.
Let $\om\colon I \to \h^\ast$ be given by
$v_i\in V_{\om(i)}$. We can write for
the right corepresentation $R\colon V\to V\wtt\A$
its matrix elements by
$Rv_i = \sum_{k\in I} v_k \otimes R_{ki}$,
$R_{ki}\in \A_{\om(k),\om(i)}$,
and the requirements $(R\otimes\id)\circ R =
(\id\otimes\De^\A)\circ R$ and
$(\id\otimes\ep^\A)\circ R=\id$ translate
into
\begin{equation}\label{eq:requirmatrixeltR}
\De^\A(R_{ij}) =\sum_{k\in I} R_{ik}\otimes R_{kj},
\qquad \ep^\A(R_{ij}) = \de_{ij} T_{-\om(i)}.
\end{equation}
Note that \eqref{eq:requirmatrixeltR} combined
with the antipode axiom
\eqref{eq:defS} implies
\begin{equation}\label{eq:SmatrixeltR}
\sum_{k\in I} S^\A(R_{ik})R_{kj} = \de_{ij} 1_\A =
\sum_{k\in I} R_{ik}S^\A(R_{kj}).
\end{equation}

Assume the existence of a form on $V$
making $R$ a unitarisable corepresentation
such that a homogeneous basis is orthogonal,
$\pair{v_i}{v_j}= \de_{ij} N_i$, $N_i\in\Mh$.
Then $\mu_r^\A(\pair{v_i}{v_j}_D\1)
=\de_{ij} \mu_r^\A(T_{\overline{\om(i)}}N_i)$,
and
\begin{equation*}
\begin{split}
\pair{Rv_i}{Rv_j} &= \sum_{k,l\in I}
\pair{v_k\otimes R_{ki}}{v_l\otimes R_{lj}} =
\sum_{k,l\in I} \pair{v_k}{v_l}_D \otimes
R_{lj}^\ast R_{ki} \\
&= \sum_{k\in I} T_{\overline{\om(k)}}\circ N_k \circ T_{-\om(k)}
\otimes R_{kj}^\ast R_{ki}
= \sum_{k\in I} R_{kj}^\ast \mu_l^\A(N_k) R_{ki}
= \de_{ij} \mu_r^\A(T_{\overline{\om(i)}}N_i)
\end{split}
\end{equation*}
using the unitarisability in the last equation.
Multiplying this last equation
by $S^\A(R_{ip})$ from the right, summing over $i\in I$
we obtain from \eqref{eq:SmatrixeltR}
\begin{equation}\label{eq:Runitaryformatrixelt}
R_{pj}^\ast \mu_l^\A(N_p) =
\mu_r^\A(T_{\overline{\om(j)}}N_j) S^\A(R_{jp})
\ \Longrightarrow\
\mu_l^\A(N_p) R_{pj} = \mu_r^\A(N_j) S^\A(R_{jp})^\ast
=S^\A( \mu_l^\A(N_j)R_{jp})^\ast.
\end{equation}
Now \eqref{eq:Runitaryformatrixelt} corresponds
to the definition of unitarisability of a matrix
corepresentation as in \cite[Def.~3.11]{KoelR}.
Note that \eqref{eq:Runitaryformatrixelt}
implies $\overline{\om(i)}=\om(i)$ for all $i\in I$,
since $S^\A(R_{jp})^\ast\in \A_{\overline{\om(p)},\overline{\om(j)}}$
and $R_{pj}\in \A_{\om(p),\om(j)}$.

Similarly, for a unitarisable left corepresentation
$Lv_i = \sum_{k\in I} L_{ik}\otimes v_k$ with
$\pair{v_i}{v_j}=\de_{ij}N_i$
we find
\begin{equation*}
\pair{Lv_i}{Lv_j} = \sum_{k\in I} L_{jk}^\ast \mu_r^\A(N_k) L_{ik}
= \de_{ij} \mu_l^\A(T_{\overline{\om(i)}} N_i).
\end{equation*}
Applying $S^\A$ to this identity, multiplying from the
left by $L_{pi}$, summing over $i\in I$, we find from
the analogue of \eqref{eq:SmatrixeltR} for the matrix elements
of the left corepresentation,
\begin{equation}\label{eq:Lunitaryformatrixelt}
S^\A\bigl( (\mu_r^\A(N_p)L_{jp})^\ast\bigr) =
\mu_l^\A( N_p) S^\A(L_{jp}^\ast) = \mu_r^\A(N_j) L_{pj}.
\end{equation}

Note that we can use the standard Gram-Schmidt process to obtain
from a set of linearly independent (over $\Mh$) vectors in $V$ a
set of orthogonal vectors, but in general not orthonormal vectors.

\begin{prop}\label{prop:actiononmatrixeltsunitarycorep}
Assume $\U$ and $\A$ are paired \hHopfalg s
having invertible antipodes. \par\noindent
{\rm (i)} Let $R$ be a
unitarisable right corepresentation of $\A$ in $V$
such that $\{ v_i\}_{i\in I}$ is a homogeneous
orthogonal basis for $V$ with $\pair{v_i}{v_j}=\de_{ij}N_j$.
Let $Rv_i = \sum_{k\in I} v_k\otimes R_{ki}$,
then
\begin{equation}\label{eq:propactiononmatrixeltsunitarycorepR}
S^\A\bigl( X^\ast \cdot (\mu_l^\A(N_j) R_{ji})\bigr)^\ast
= (\mu_l^\A(N_i) R_{ij})\cdot X, \qquad \forall X\in \U .
\end{equation}
\par\noindent
{\rm (ii)} Let $L$ be a
unitarisable left corepresentation of $\A$ in $V$
such that $\{ v_i\}_{i\in I}$ is a homogeneous
orthogonal basis for $V$ with $\pair{v_i}{v_j}=\de_{ij}N_j$.
Let $Lv_i = \sum_{k\in I} L_{ik}\otimes v_k$,
then
\begin{equation}\label{eq:propactiononmatrixeltsunitarycorepL}
 X \cdot (\mu_r^\A(N_j) L_{ij})
= S^\A\bigl( ((\mu_r^\A(N_i) L_{ji})\cdot X^\da)^\ast\bigr),
\qquad \forall X\in \U .
\end{equation}
\end{prop}

Note that \eqref{eq:propactiononmatrixeltsunitarycorepR},
respectively \eqref{eq:propactiononmatrixeltsunitarycorepL},
reduces to \eqref{eq:Runitaryformatrixelt},
respectively \eqref{eq:Lunitaryformatrixelt},
for $X=1\in \U$.

\begin{proof}  Start with
$R_{ij}\cdot X = \sum_{k\in I} T_\al\pair{X}{R_{ik}}\otimes R_{kj}$,
$X\in\U_{\al\be}$,
in which we use \eqref{eq:Runitaryformatrixelt},
\eqref{eq:pairing-leftmoment},
\eqref{eq:pairing-rightmoment} and
\eqref{eq:pairingXstara} to find the first equality in
\begin{equation*}
\begin{split}
R_{ij}\cdot X &=
\sum_{k\in I} T_\al\circ  N_i^{-1}\circ T_{-\om(i)}
\circ \pair{X^\ast}{R_{ki}}^\ast \circ N_k \circ T_{-\om(k)}
\otimes \mu_r^\A(N_j)\mu_l^\A(N_k^{-1}) S^\A(R_{jk})^\ast \\
& =
\Bigl( \sum_{k\in I} T_{\om(k)}\circ  N_k\circ
\pair{X^\ast}{R_{ki}} \circ T_{\om(i)}\circ
N_i^{-1} \circ T_{-\bar\al}
\otimes S^\A(R_{jk})\mu_r^\A(N_j)\mu_l^\A(N_k^{-1}) \Bigr)^\ast \\
& = S^\A \Bigl( \sum_{k\in I}
\mu_l^\A(N_j)\mu_r^\A(N_k^{-1})R_{jk} \otimes
T_{\bar\al}\circ  N_i^{-1}\circ T_{-\om(i)}
S^\Dh(\pair{X^\ast}{R_{ki}}) \circ N_k\circ T_{-\om(k)}\Bigr)^\ast.
\end{split}
\end{equation*}
Write $\pair{X^\ast}{R_{ki}}=fT_{\bar\al-\om(i)}$, so
$S^\Dh(\pair{X^\ast}{R_{ki}})=T_{-\bar\al+\om(i)} \circ f$,
to find using \eqref{eq:identDhotimesAwithA}
\begin{equation*}
\begin{split}
R_{ij}\cdot X &= S^\A\Bigl( \sum_{k\in I}
\mu_r^\A \bigl( (T_{\bar\al}N_i^{-1}) f N_k )
\mu_l^\A(N_j)\mu_r^\A(N_k^{-1})R_{jk}\Bigr)^\ast \\
&= \mu_l^\A (T_{\al}N_i^{-1}) \mu_r^\A(N_j)
\ S^\A\Bigl( \sum_{k\in I}
\mu_r^\A( \pair{X^\ast}{R_{ki}}T_{-\bar\be}\1) R_{jk}\Bigr)^\ast,
\end{split}
\end{equation*}
which we rewrite as
$\mu_l^\A (T_{\al}N_i) \bigl( R_{ij}\cdot X\bigr)
= S^\A\bigl( \mu_l^\A(N_j) (X^\ast \cdot R_{ji})\bigr)^\ast$
and the result follows from \eqref{eq:actionmomentA}.

The statement (ii) for the left corepresentation is
proved analogously.
\end{proof}

\begin{prop}\label{prop:repforunitarycorep}
Assume $\U$ and $\A$ are paired \hHopfster s having
invertible antipodes. \par\noindent {\rm (i)} Let $R$ be a
unitarisable right corepresentation of $\A$ in $V$, and let
$\pi(X)$, $X\in\U_{\al\be}$, be defined as in Proposition
\ref{prop:cotodynrep}(i), then $T_\be\circ \pair{\pi(X)v}{w}_D =
\pair{v}{\pi(X^*)w}_D\circ T_\be$.
\par\noindent
{\rm (ii)} Let $L$ be a unitarisable left corepresentation of $\A$
in $V$, and let $\pi(X)$, $X\in\U_{\al\be}$, be defined as in
Proposition \ref{prop:cotodynrep}(ii), then
$\pair{\pi(X)v}{w}_D\circ T_\be = T_\be\circ
\pair{v}{\pi(X^\da)w}_D$.
\end{prop}

\begin{proof} For (i) note that both sides are not sesquilinear
over $\Mh$, but for $f\in\Mh$ and $v\in V_\ga$ we have
$\pair{fv}{\pi(X^*)w}_D\circ T_\be = \pair{v}{\pi(X^*)w}_D\circ
(T_\ga f) \circ T_\be = \pair{v}{\pi(X^*)w}_D \circ T_\be \circ
(T_{\ga-\be} f)$, and on the other hand $T_\be\circ
\pair{\pi(X)(fv)}{w}_D = T_\be\circ
\pair{(T_{-\al}f)\pi(X)v}{w}_D= T_\be\circ
\pair{\pi(X)v}{w}_D\circ (T_{\ga-\be}f)$ using
$\pi(X)V_\ga\subseteq V_{\ga+\al-\be}$. Using $\pair{v}{w}_D
=\pair{w}{v}_D^\ast$, or performing a similar calculation we find
$T_\be\circ \pair{\pi(X)v}{f w}_D = (T_{\bar\de+\be}\bar f) \circ
T_\be \circ \pair{\pi(X)v}{w}_D$ and
$\pair{v}{\pi(X^*)(fw)}_D\circ T_\be = (T_{\bar\de+\be}\bar f)
\circ \pair{v}{\pi(X^*)w}_D\circ T_\be$. This shows that we can
restrict to proving the result for $v$, $w$ elements of a basis of
$V$ over $\Mh$.

Let $\{v_i\}_{i\in I}$, $v_i\in V_{\om(i)}$,
be a homogeneous basis for $V$ over $\Mh$
with $\pair{v_i}{v_j}=\de_{ij}N_j$,
$Rv_i=\sum_{k\in I} R_{ik}\otimes v_k$. Then the
unitarisability of $R$ is expressed by
\eqref{eq:Runitaryformatrixelt}. Now
\begin{equation*}
\begin{split}
\pair{v_i}{\pi(X^*)v_j}_D &= \pair{v_i\otimes
T_{-\om(i)}}{\sum_{l\in I} v_l\otimes \pair{X^*}{R_{lj}}\circ
T_{-\bar\be}} = \sum_{l\in I} \pair{v_i}{v_l}_D \otimes
T_\be \circ \pair{X^*}{R_{lj}}^\ast \circ T_{-\om(i)} \\
&= \sum_{l\in I} \pair{v_i}{v_l}_D \otimes
T_\be \circ T_{\overline{\om(j)}}\circ
\pair{X}{S^\A(R_{lj})^\ast}\circ T_{\overline{\om(l)}} \circ T_{-\om(i)}
\end{split}
\end{equation*}
and using \eqref{eq:Runitaryformatrixelt} and
\eqref{eq:pairing-leftmoment}, \eqref{eq:pairing-rightmoment}, and
the orthogonality of the basis we find
\begin{equation*}
\begin{split}
\pair{v_i}{\pi(X^*)v_j}_D &= T_{\overline{\om(i)}}\circ N_i \circ
T_{-\om(i)} \otimes T_\be \circ T_{\om(j)} \circ N_j \circ
\pair{X}{R_{ji}} \circ (T_{\om(i)}N_i^{-1}) \circ
T_{\overline{\om(i)}}\circ T_{-\om(i)} \\
&= (T_{\overline{\om(i)}} N_i) \circ T_{\be+\om(j)} \circ N_j
\circ \pair{X}{R_{ji}} \circ (T_{\om(i)}N_i^{-1})
= T_{\be+\om(j)} \circ N_j
\circ \pair{X}{R_{ji}}
\end{split}
\end{equation*}
using $\overline{\om(i)}=\om(i)$ and $T_{\be+\om(j)}\circ N_j
\circ \pair{X}{R_{ji}}$ being a multiplication operator
as well.
Comparing this with
\begin{equation*}
\begin{split}
\pair{\pi(X)v_i}{v_j} &= \sum_{k\in I} \pair{v_k}{v_j}_D \otimes
T_{\overline{\om(j)}}\circ \pair{X}{R_{ki}}\circ T_\be
= T_{\overline{\om(j)}}N_j T_{-\om(j)} \otimes
T_{\overline{\om(j)}}\circ \pair{X}{R_{ji}}\circ T_\be \\
&= (T_{\overline{\om(j)}}N_j) \circ
T_{\overline{\om(j)}}\circ \pair{X}{R_{ji}}\circ T_\be
= T_{\om(j)}\circ N_j \circ \pair{X}{R_{ji}}\circ T_\be
\end{split}
\end{equation*}
proves the first statement. Part (ii) is proved along the same lines.
%
\end{proof}

\subsection{Elliptic $U(2)$ dynamical quantum group}
\label{ssec:ellpairingmatrixelt}

In this subsection we calculate the pairing between matrix
elements of irreducible corepresentations of the
elliptic $U(2)$ dynamical quantum group using the pairing
on the generators given by the $R$-matrix of
\eqref{eq:ellR} as in  Remark \ref{rmk:cobraidingandRmatrix}(ii).
It turns out that the pairing is described in terms of
elliptic hypergeometric series, whose definition we recall
later. Such series have been considered for the first time
by Frenkel and Turaev \cite{FrenT}, and
this paper has triggered a lot of activity, see
\cite[Ch. 11]{GaspR}
and references given there. We show how the link between the
pairing of matrix elements can be used to rederive the results
of Frenkel and Turaev \cite{FrenT}, giving a somewhat alternative
quantum group theoretic derivation
of these results for elliptic hypergeometric series
as in  \cite{KoelvNR}.

Before doing so we recall the corepresentations of the
\hHopfster\ $\E$ as constructed in \cite{KoelvNR}.
Let $N\in\N$ and put
\begin{equation}\label{eq:defvkz}
v_k(z)= \gamma(z)\gamma(q^2z)\cdots\gamma(q^{2(N-k-1)}z)
\alpha(q^{2(N-k)}z)\cdots\alpha(q^{2(N-1)}z), \qquad
k\in\{0,1,\ldots,N\},
\end{equation}
and let $V^N_{2k-N}=\mu^\E_l(M_{\h^*})v_k(z)$,
$V^N=\bigoplus_{k=0}^N V^N_{2k-N}$. Then $V^N$ is an
\hspa\  with  $\mu_{V^N}$
given by the left moment map $\mu^\E_l$.
Then $V^N$ is a left corepresentation of $\E$ given
by the restriction of the comultiplication;
\begin{equation}\label{eq:defleftcorepU}
L \colon V^N \to \E\otimes V^N, \qquad L= \De^\E\big\vert_{V^N}, \qquad
\De^\E(v_k(z)) = \sum_{i=0}^N t^N_{kj}(z) \otimes v_j(z).
\end{equation}
The matrix elements can be expressed as follows,
see \cite[Thm. 3.4]{KoelvNR}
for the first expression and its proof for the second expression;
\begin{equation}\label{eq:explexprtNkj}
\begin{split}
t^N_{kj}(z)=&\sum_{l=\max(0,k+j-N)}^{\min(k,j)}
\begin{bmatrix} k\\l \end{bmatrix}\begin{bmatrix} N-k\\j-l \end{bmatrix}
\mu_r^\E \Bigl( \frac{(q^{2(\cdot +N-k-2j+l+2)})_l}{(q^{2(\cdot +N-2j+2)})_l}
 \frac{(q^{2(\cdot +l-j+2)})_{j-l}}{(q^{2(\cdot +N-2j-k+2l+2)})_{j-l}}\Bigr)\\
    &\qquad\times\gamma(q^{2(N-k-1)}z)\cdots\gamma(q^{2(N-j-k+l)}z)
\delta(q^{2(N-j-k+l-1)}z)\cdots\delta(z)\\
    &\qquad\times\alpha(q^{2(N-1)}z)\cdots\alpha(q^{2(N-l)}z)
\beta(q^{2(N-l-1)}z)\cdots\beta(q^{2(N-k)}z) \\
=& \sum_{l=\max(0,j+k-N)}^{\min(k,j)}
t^{N-k}_{0,j-l}(z)\, t^k_{kl}(q^{2(N-k)}z),
\end{split}
\end{equation}
where the elliptic binomial coefficient is defined by
$\displaystyle{\begin{bmatrix} k\\l \end{bmatrix}=\prod_{i=1}^l
\frac{\theta(q^{2(k-l+i)})}{\theta(q^{2i})}}$,
and the elliptic shifted factorials are
\begin{equation}\label{eq:ellshiftedfactorial}
(a)_n = \prod_{i=0}^{n-1} \theta(aq^{2i}), \qquad
(a_1,\ldots, a_k)_n =\prod_{i=1}^k (a_i)_n.
\end{equation}
It then follows that $t^N_{kj}(z)\in \E_{2k-N,2j-N}$,
$\ep^\E(t^N_{kj}(z))=\de_{kj}T_{N-2k}$, $\De^\E(t^N_{kj}(z))=\sum_{p=0}^N
t^N_{kp}(z)\otimes t^N_{pj}(z)$.

We first need to calculate the pairing of the generators with the
matrix elements.

\begin{thm}\label{thm:pairingellgenmatrixelt}
The cobraiding $\E^\cop\times\E\to\Dh$
satisfies
\begin{equation*}
\begin{split}
\pair{\al(w)}{t^N_{kj}(z)} &=
\de_{k,j} \frac{\theta(q^{2(1-N+k)}w/z, q^{2(\la+N-k+1)})}
{\theta(q^2w/z, q^{2(\la+1)})}\, T_{N-2k-1}, \\
\pair{\be(w)}{t^N_{kj}(z)} &=
\de_{k,j-1}\frac{\theta(q^{2(N-k)}, q^{-2(\la+N-k)}w/z)}
{\theta(q^2w/z, q^{-2(\la+1)})}\, T_{N-2k-1}, \\
\pair{\ga(w)}{t^N_{kj}(z)} &=
\de_{j,k-1} \frac{
\theta(q^{2k}, q^{2(\la-k+2)}w/z)}
{\theta(q^2w/z, q^{2(\la+1)})} \, T_{N-2k+1}, \\
\pair{\de(w)}{t^N_{kj}(z)} &=
\de_{k,j} \frac{\theta(q^{2(k-\la-1)},q^{2(1-k)}w/z)}
{\theta(q^2w/z, q^{-2(\la+1)})} T_{N-2k+1}, \\
\pair{\det^{-1}(w)}{t^N_{kj}(z)} & =\de_{k,j}\ q^{-N}
\frac{\theta(q^2w/z)}{\theta(q^{2(1-N)}w/z)}\, T_{N-2k},
\end{split}
\end{equation*}
where we consider the functions in $\la$ on the right hand
sides as elements of $\Mh$.
\end{thm}

\begin{proof}
First observe that $t^N_{kj}(z) \in (\E)_{2k-N,2j-N}$, and
$\al(w)\in \E_{11}=(\E^\cop)_{11}$, $\be(w)\in
\E_{1,-1}=(\E^\cop)_{-1,1}$, $\ga(w)\in
\E_{-1,1}=(\E^\cop)_{1,-1}$, and $\de(w)\in
\E_{-1,-1}=(\E^\cop)_{-1,-1}$. So $(\Dh)_{\al\be}=\{0\}$
for $\al\not=\be$ and \eqref{eq:pairing-target} imply
$\pair{\al(w)}{t^N_{kj}(z)}=0$ unless $1+2j-N=1+2k-N$ or
$k=j$, similarly $\pair{\be(w)}{t^N_{kj}(z)}=0$ unless
$j-k=1$, $\pair{\ga(w)}{t^N_{kj}(z)}=0$ unless $k-j=1$, and
$\pair{\de(w)}{t^N_{kj}(z)}=0$ unless $k=j$.

From the explicit expression \eqref{eq:explexprtNkj}
we see that $t^N_{kj}(z)$ is
composed of elements of the form ($k\geq 1$)
\begin{equation}\label{eq:matcorepcorners}
\begin{split}
t^k_{00}(z) &= \de(q^{2(k-1)}z)\cdots \de(z), \qquad
t^k_{0k}(z) = \ga(q^{2(k-1)}z)\cdots \ga(z), \\
t^k_{k0}(z) &= \be(q^{2(k-1)}z)\cdots \be(z), \qquad
t^k_{kk}(z) = \al(q^{2(k-1)}z)\cdots \al(z).
\end{split}
\end{equation}
Using \eqref{eq:pairing-rightproduct}, the comultiplication for
the generators defined by \eqref{eq:epandDeFRST} (see
\cite{FeldV}, \cite{KoelvNR} for explicit formulas),
\eqref{eq:defabcd},
\eqref{eq:elldualgenerator} we prove by induction on $k$
that the only non-zero pairings between generators and
matrix elements of the form \eqref{eq:matcorepcorners}
are given by
\begin{equation}\label{eq:pairinggenmatcorepcorners}
\begin{split}
\pair{\al(w)}{t^k_{00}(z)} &=
\frac{\theta(q^{2(1-k)}w/z, q^{2(\la+k+1)})}
{\theta(q^2w/z, q^{2(\la+1)})}\, T_{k-1}, \qquad
\pair{\al(w)}{t^k_{kk}(z)} = T_{-k-1}, \\
\pair{\de(w)}{t^k_{kk}(z)} &=
\frac{\theta(q^{-2(k-1)}w/z, q^{-2(\la-k+1)})}
{\theta(q^2w/z, q^{-2(\la+1)})}\, T_{1-k},
 \qquad
\pair{\de(w)}{t^k_{00}(z)} = T_{k+1}, \\
\pair{\be(w)}{t^k_{0k}(z)} &=
\begin{cases} b(\la,w/z)\, T_0, & k=1,\\
0, & k>1, \end{cases}, \qquad
\pair{\ga(w)}{t^k_{k0}(z)} =
\begin{cases} c(\la,w/z)\, T_0, & k=1,\\
0, &k>1. \end{cases}
\end{split}
\end{equation}
This proves the theorem in this particular case.

Next we treat the pairing with $t^N_{Nj}(z)$.
Now $\pair{\be(w)}{t^N_{Nj}(z)}=0$ for all $j$
by the weight considerations in the first paragraph
of the proof, and the non-zero pairings with $\al(w)$,
$\de(w)$ are already contained
in \eqref{eq:pairinggenmatcorepcorners} for $j=N$. The
only non-zero case is $\pair{\ga(w)}{t^N_{N,N-1}(z)}$,
and from \eqref{eq:explexprtNkj} we see that
$t^N_{N,N-1}(z)$ is a product of $\mu_r^\E(C_N)$ for
an explicit function $C_N$ times $N-1$ $\al$'s and one
$\be$. Bringing the function to the right, using
\eqref{eq:pairing-rightmoment}, and next using the comultiplication
$\De^{\E^\cop}(\ga(w))=\al(w)\otimes\ga(w) + \ga(w)\otimes
\de(w)$, we calculate the pairing from
\eqref{eq:pairinggenmatcorepcorners} and
\eqref{eq:elldualgenerator}, where we write the functions
in $\Mh$ as functions of $\la$;
\begin{equation*}
\begin{split}
\pair{\ga(w)}{t^N_{N,N-1}(z)} &=
\pair{\ga(w)}{\al(q^{2(N-1)}z)\cdots
\al(q^2z)\be(z)} \circ (T_{N-2}C_N)(\la) \\
&= \Bigl( \pair{\al(w)}{\al(q^{2(N-1)}z)\cdots \al(q^2z)} T_1
\pair{\ga(w)}{\be(z)} \\
&\qquad + \pair{\ga(w)}{\al(q^{2(N-1)}z)\cdots \al(q^2z)}
T_{-1} \pair{\de(w)}{\be(z)}\Bigr) \circ C_N(\la+N-2)  \\
&= T_{-(N-1)-1}\circ T_1 \circ c(\la,\frac{w}{z}s) C_N(\la+N-2) =
\frac{\theta(q^{2N}, q^{2(\la-N+2)}w/z)}
{\theta(q^2w/z, q^{2(\la+1)})}
\, T_{-N+1}.
\end{split}
\end{equation*}
Similarly, we can establish the pairing of a generator
with $t^N_{0j}(z)$. Now
$\pair{\ga(w)}{t^N_{0j}(z)}=0$ for all $j$, and the
non-zero pairings with $\al(w)$, $\de(w)$ are already contained
in \eqref{eq:pairinggenmatcorepcorners} for $j=0$. The
only non-zero case is $\pair{\be(w)}{t^N_{01}(z)}$, and
we prove similarly
\begin{equation*}
\pair{\be(w)}{t^N_{01}(z)} = \frac{\theta(q^{2N}, q^{-2(\la+N)}w/z)}
{\theta(q^2w/z, q^{-2(\la+1)})}\, T_{N-1}.
\end{equation*}

After these preparations we can use the second expression
of \eqref{eq:explexprtNkj} to calculate the pairing of a
generator with an arbitrary matrix element. We give the
details for the pairing $\pair{\de(w)}{t^N_{kj}(z)}$.
\begin{equation*}
\begin{split}
\pair{\de(w)}{t^N_{kj}(z)} &=
\sum_{l=\max(0,j+k-N)}^{\min(k,j)}
\pair{\be(w)}{t^{N-k}_{0,j-l}(z)} \,  T_1 \,
\pair{\ga(w)}{t^k_{kl}(q^{2(N-k)}z)} \\
&\ \ \ +
\sum_{l=\max(0,j+k-N)}^{\min(k,j)}
\pair{\de(w)}{t^{N-k}_{0,j-l}(z)} \, T_{-1}\,
\pair{\de(w)}{t^k_{kl}(q^{2(N-k)}z)}.
\end{split}
\end{equation*}
The summand in the first sum is non-zero if and
only if $j-l=1$ and $l=k-1$, and the summand in the
second sum is non-zero if and only if
$j=l$ and $k=l$, so we are left with only two non-zero
terms in the case $k=j$ and no non-zero terms in
case $k\not=j$.
For the pairing with $\al(w)$, $\be(w)$
and $\ga(w)$ instead of $\de(w)$ we only get at most
one non-zero term, so that these case are simpler.
So we obtain $\pair{\de(w)}{t^N_{kj}(z)}$
\begin{equation*}
\begin{split}
=& \de_{kj}
\Bigl( \pair{\be(w)}{t^{N-k}_{01}(z)} \,  T_1 \,
\pair{\ga(w)}{t^k_{k,k-1}(q^{2(N-k)}z)} 
+ \pair{\de(w)}{t^{N-k}_{00}(z)} \, T_{-1}\,
\pair{\de(w)}{t^k_{kk}(q^{2(N-k)}z)}\Bigr) \\
=& \de_{kj}\Bigl(
\frac{\theta(q^{2(N-k)}, q^{-2(\la+N-k)}w/z)}
{\theta(q^2w/z, q^{-2(\la+1)})}\, T_{N-k-1}\, T_1\,
\frac{\theta(q^{2k}, q^{2(\la-N+2)}w/z)}
{\theta(q^{2(1+k-N)}w/z, q^{2(\la+1)})}
\, T_{-k+1} \\
&\qquad\qquad\qquad\qquad\qquad\qquad\qquad\qquad
+ T_{N-k+1}\, T_{-1}\,
\frac{\theta(q^{2(1-N)}w/z, q^{-2(\la-k+1)})}
{\theta(q^{2(1+k-N)}w/z, q^{-2(\la+1)})}\, T_{1-k}\Bigr)
\\
=& \de_{kj}\Bigl(
\frac{\theta(q^{2(N-k)}, q^{-2(\la+N-k)}w/z,q^{2k}, q^{2(\la-k+2)}w/z)}
{\theta(q^2w/z, q^{-2(\la+1)},q^{2(1+k-N)}w/z, q^{2(\la+N-k+1)})}
+ \frac{\theta(q^{2(1-N)}w/z, q^{-2(\la+N-2k+1)})}
{\theta(q^{2(1+k-N)}w/z, q^{-2(\la+N-k+1)})}
\Bigr)\, T_{N-2k+1}.
\end{split}
\end{equation*}
It remains to calculate the term in parentheses, and this is
done using the addition theorem \eqref{eq:addell} with
$(x,y,z,w)$ replaced by $(q^{-\la}\sqrt{w/z},
q^{2(N-k)+\la}\sqrt{z/w}, q^{\la+2}\sqrt{w/z},
q^{2k-\la-2}\sqrt{z/w})$ to rewrite the numerator
of the first quotient. A straightforward calculation
gives the required result. We note that the calculation
for the pairing with $\al(w)$, $\be(w)$
and $\ga(w)$ instead of $\de(w)$ does not require
the addition formula \eqref{eq:addell}.

Having the pairings with the generators $\al(w)$,
$\be(w)$, $\ga(w)$, $\de(w)$, we can calculate the
pairing with $\det(w)$ and from this derive the
pairing with $\det^{-1}(w)$ as in \S \ref{ssec:ellpairing}
using \eqref{eq:addell} again.
\end{proof}

Now we combine Theorem \ref{thm:pairingellgenmatrixelt}
with Proposition \ref{prop:cotodynrep}
for the corresponding left corepresentation $V^N$
as in \eqref{eq:defvkz}, \eqref{eq:defleftcorepU}.
This yields the
dynamical representation of $\E^{\opp}$ on $V^N$
given by
\begin{equation}\label{eq:repAcopopp}
\begin{split}
\pi(\al(w)) \, \bigl( \mu_{V^N}(f)v_k(z)\bigr) &=
\mu_{V_N}\left(\frac{\theta(q^{2(1-N+k)}w/z, q^{2(\la+N-k+2)})}
{\theta(q^2w/z, q^{2(\la+2)})}\, (T_1f)\right)\,  v_k(z), \\
\pi(\be(w))  \, \bigl( \mu_{V^N}(f)v_k(z)\bigr) &=  \mu_{V_N}\left(\frac{\theta(q^{2(N-k)},
q^{-2(\la-1+N-k)}w/z)}
{\theta(q^2w/z, q^{-2\la})} \, (T_{-1}f)\right) \, v_{k+1}(z), \\
\pi(\ga(w))  \, \bigl( \mu_{V^N}(f)v_k(z)\bigr) &=
\mu_{V_N}\left(\frac{
\theta(q^{2k}, q^{2(\la-k+3)}w/z)}
{\theta(q^2w/z, q^{2(\la+2)})}
\, (T_1f)\right)\, v_{k-1}(z), \\
\pi(\de(w))  \, \bigl( \mu_{V^N}(f)v_k(z)\bigr) &=
\mu_{V_N}\left( \frac{\theta(q^{2(k-\la)},q^{2(1-k)}w/z)}
{\theta(q^2w/z, q^{-2\la})}\, (T_{-1}f)\right) \,  v_k(z), \\
\pi(\det^{-1}(w)) \, \bigl( \mu_{V^N}(f)v_k(z)\bigr)
&= q^{-N} \frac{\theta(q^2w/z)}{\theta(q^{2(1-N)}w/z)}\, \mu_{V_N}(f) v_k(z).
\end{split}
\end{equation}

We now proceed as follows to determine the
pairing of matrix elements. First we use the
dynamical representation of $\E^\opp$
in \eqref{eq:repAcopopp} and the explicit
expression \eqref{eq:explexprtNkj} to
calculate $\pi\big( t^M_{rs}(w)\bigr)\, v_k(z)$
explicitly as a multiple of $v_{k-s+r}(z)$,
whereas from the definition in Proposition
\ref{prop:cotodynrep} we have
\begin{equation*}
\pi(t^M_{rs}(w))\, v_k(z) =
\sum_{j=0}^N T_{2s-M} \pair{t^M_{rs}(w)}{t^N_{kj}(z)}
\otimes v_j(z) =
\sum_{j=0}^N \mu_V\bigl( T_{2s-M}
\pair{t^M_{rs}(w)}{t^N_{kj}(z)}\1\bigr) \, v_j(z),
\end{equation*}
and upon comparing these expression we obtain the
desired result in Theorem \ref{thm:pairingmatrixeltE}.
In order to state the result we need some notations from
special functions in a special case. The
(very-well-poised) elliptic
hypergeometric series is defined by
\begin{equation}\label{eq:ellhypseries}
{}_{r+1}V_r( a_1;a_6,\ldots, a_{r+1}) =
\sum_{n=0}^\infty \frac{\theta(a_1q^{4n})}{\theta(a_1)}
\frac{ q^{2n} (a_1,a_6,a_7,\ldots,a_{r+1})_n}
{(q^2,q^2a_1/a_6,q^2a_1/a_6,\ldots, q^2a_{r+1}/a_6)_n},
\end{equation}
assuming the elliptic balanced condition $(a_6 a_7\ldots a_{r+1})^2 q^4 = (a_1q^2)^{r-5}$ and
where we have used the notation \eqref{eq:ellshiftedfactorial}.
This is not the most general definition, but it
suffices for our purposes, see \cite[\S 11.2]{GaspR}, \cite{Spir}
for a discussion and references. Compared to these references
we have switched from base $q$ to base $q^2$, and we have
specialised $z=1$. In this paper all elliptic hypergeometric
series are terminating series.

\begin{thm}\label{thm:pairingmatrixeltE}
The pairing between matrix elements of the irreducible
corepresentations as in \eqref{eq:explexprtNkj} are given by
\begin{equation}\label{eq:pairingmatrixeltE}
\begin{split}
&\pair{t^M_{rs}(w)}{t^N_{kj}(z)}
= \de_{s+j,r+k} C\
{}_{12}V_{11}\bigl(q^{2(\la+M-2s-r+1)};
 q^{-2r},q^{-2s}, q^{2(\la-s+1)},
  q^{2(\la+M+N-k-s-r+2)}, \\
  & \qquad\qquad\qquad\qquad\qquad \qquad\qquad q^{2(\la+M-k-s-r+1)},
q^{2(k-s)}z/w, q^{2(M-N+k-s+1)}w/z\bigr) \,
T_{N+M-2s-2j}, \\
&C=\;(-1)^{M-r-s}q^{(M-r-s)(\la-s+1)}
\frac{(q^{2(M-r-s+1)},q^{2(k-s+1)},q^{2(\la+M-s-k-r+2)}w/z)_s}
{(q^2,q^{2(\la+M-2s-r+2)},q^{2(M-r-s+1)}w/z)_s}\\
 & \qquad \times\frac{(q^{2(N-k+s-r+1)},q^{2(-\la-N+k+s)}w/z)_r}{(q^{-2(\la+M-2s)},
 q^{2(M-r+1)}w/z)_r}
\frac{(q^{2(k+s+r-\la-M)},q^{2(s-k+1)}w/z)_{M-r-s}}{(q^{2(\la-s+1)},q^2
w/z )_{M-r-s}}.
\end{split}
\end{equation}
\end{thm}

Before going into the proof of Theorem
\ref{thm:pairingmatrixeltE}, we show how the
quantum dynamical Yang-Baxter equation, the
biorthogonality relations, Bailey's transformation
formula, and
it is well known that a corollary to Bailey's transformation is
Jackson's summation formula, are implied by
Theorem \ref{thm:pairingmatrixeltE}. Since we re-obtain
the well-known properties of the elliptic hypergeometric
series already obtained by Frenkel and Turaev \cite{FrenT},
we only sketch the derivation.

First, by taking $a=t^M_{rs}(w)$, $b=t^N_{kj}(z)$ in
the cobraiding property \eqref{eq:defcobraiding} we
obtain an identity in $\E$, pairing this identity
with an element $t^L_{ac}(u) \in \E^\cop$ gives an
identity in $\Dh$.
By Theorem \ref{thm:pairingmatrixeltE}
the resulting identity is trivial unless
$r+k+a=s+j+b$. In the latter case, it gives an
identity of the form
$\sum_{\text{single}} {}_{12}V_{11}\, {}_{12}V_{11}\, {}_{12}V_{11} =
\sum_{\text{single}} {}_{12}V_{11}\, {}_{12}V_{11}\, {}_{12}V_{11}$,
which can be rephrased as an $R$-matrix satisfying
the quantum dynamical Yang-Baxter
equation \eqref{eq:QDYBE}, the case $L=M=N=1$ corresponding to the
$R$-matrix \eqref{eq:ellR}, see \cite{FrenT}, or as an
elliptic analogue of Wigner's symmetry for the
$9j$-symbols.

The proof of the biorthogonality relations for the elliptic
hypergeometric series is essentially the same as in
\cite{KoelvNR}, i.e. we pair $\sum_{k=0}^N t^N_{jk}(z)\, S^\A(t^N_{ki}(z)) =\de_{ij}1_{\E}$ with an arbitrary $t^M_{rs}(w)$. For this we
also need the pairing $\pair{t^M_{rs}(w)}{S^\A(t^N_{kj}(z))}$, which
can be calculated using the unitarity of the corepresentation
$t^N$ as in \cite{KoelvNR}. Also, Bailey's transform can be
obtained from the unitarity of the corepresentation $t^N$, see
\cite{KoelvNR}.

\begin{proof}[Proof of Theorem \ref{thm:pairingmatrixeltE}]
Observe that $T_{2s-M} \pair{t^M_{rs}(w)}{t^N_{kj}(z)}
\in \bigl(\Dh\bigr)_{2j-N,2r+2k-N-2s}$ so that there
is only a non-zero contribution for $j+s=k+r$. As noted
previously, it suffices to calculate
$\pi\big( t^M_{rs}(w)\bigr)\, v_k(z)$
explicitly as a multiple of $v_{k-s+r}(z)$.
Now By Proposition
\ref{prop:cotodynrep}(ii) $\pi$ is an antimultiplicative
representation of $\E^\cop$, so that
\eqref{eq:explexprtNkj} gives
$\pi(t^M_{rs}(w))\, v_k(z)$
\begin{equation}\label{eq:pairinghelp1}
\begin{split}
=&\sum_{l=\max(0,r+s-M)}^{\min(r,s)}
\pi(\be(q^{2(M-r)}w)) \cdots \pi(\be(q^{2(M-l-1)}w))
\pi(\al(q^{2(M-l)}w))\cdots
\pi(\al(q^{2(M-1)}w)) \\
&\qquad\times\pi(\de(w))\cdots \pi(\de(q^{2(M-s-r+l-1)}w))
\pi(\ga(q^{2(M-s-r+l)}w)) \cdots
\pi(\ga(q^{2(M-r-1)}w)) \\
&\qquad\times \pi(
\mu_l^\E\Bigl( \frac{(q^{2(\cdot +M-r-2s+l+2)})_l}
{(q^{2(\cdot +M-2r+2)})_l} \frac{(q^{2(\cdot+l-s+2)})_{s-l}}
{(q^{2(\cdot+M-2s-r+2l+2)})_{s-l}}
\begin{bmatrix} M-r \\ s-l\end{bmatrix}
\begin{bmatrix} r \\ l\end{bmatrix}\Bigr))\ v_k(z).
\end{split}
\end{equation}
It is now a tedious, but straightforward verification
using \eqref{eq:repAcopopp} that we can rewrite
\eqref{eq:pairinghelp1} as
an elliptic hypergeometric series of the type
\eqref{eq:ellhypseries} times $v_{k-s+r}$. Since
\eqref{eq:pairinghelp1} is also equal to
$\sum_{j=0}^N T_{2s-M} \pair{t^M_{rs}(w)}{t^N_{kj}(z)}
\otimes v_j(z)$, we find the result.
\end{proof}

\section{Singular and spherical vectors}\label{sec:singular}

From the previous sections it is clear that we can
calculate all pairings for the case of the elliptic
$U(2)$ dynamical quantum group in detail. In order to
deal with more general cases we want to have the
dynamical analogue of notions as spherical functions
on symmetric spaces which we want to realize as
subalgebras satisfying certain invariance conditions
stated in terms of the actions defined in
Theorem \ref{thm:defaction}.
The purpose of this section is to start a
description for these notions in the setting of dynamical
quantum groups, and to give the details for
the case of the elliptic $U(2)$ dynamical quantum group.
We expect to deal with more involved examples in the
future for which the notions in \S \ref{ssec:hsubalgebras}
are required.

\subsection{Subalgebras of \halg s}\label{ssec:hsubalgebras}

Using the notions of \S \ref{ssec:algebraicnotions}
we say that the \hprealg\ $\B$ is \emph{$\h$-presubalgebra}
of the \halg\ $\A$ if there exists injective
\hprealg\ homomorphism $\io\colon \B\to \A$, and
similarly a \halg\ $\B$ is a \emph{$\h$-subalgebra} of the
\halg\ $\A$ if there exists an injective unital
\halg\ homomorphism $\io\colon \B\to \A$.

For later use we need the notion of subalgebras of
\hprealg s and \halg s
with not necessarily the same Lie algebra $\h$.
We do not need this for the particular example
of the elliptic $U(2)$ dynamical quantum group, but
we introduce the notion for use in future work.
Assume $\rl$ is a complex vector space that
is a subspace of the complex vector space
$\h$. Let $i\colon \rl\to\h$ be the
corresponding injection.
Then the restriction $r\colon \hs \to \rls$ is given
by $r(\al)(X)=\al(i(X))$ for $\al\in\hs$,
$X\in\rl$. Note $r$ is a surjective linear map.
For a function $f$ on $\rls$
we define the function $j(f)$ on $\hs$
by $j(f)(\al)=f(r(\al))$.
It follows that $j\colon \Mr\to\Mh$.
We assume this situation in the remainder
of this section.

For a \hprealg\ $\A$ we define the {\em weights}
as $w(\A) = \{ \al\in\hs\mid \exists\, \be\in\hs
\text{ such that }\A_{\al\be}\not=\{0\} \text{ or }
\A_{\be\al}\not=\{ 0\} \}\subset \hs$.

A $\rl$-prealgebra $\B$ is a \emph{$(\rl,\h)$-presubalgebra} of
the \hprealg\ $\A$ if $\rl$ is a subspace of $\h$, and
such that there exists a
map $s\colon w(\B) \to w(\A)$ with $r\circ s=\id_{w(\B)}$
and an injective $\C$-linear map
$\io \colon \B\to \A$
such that
$\io(\B_{\si,\tau})\subseteq \A_{s(\si),s(\tau)}$
for all $\si,\tau\in w(\B)$, and such that
$\io(\mu_l^\B(f)b) = \mu_l^\A(j(f))\io(b)$ and
$\io(\mu_r^\B(f)b) = \mu_r^\A(j(f))\io(b)$
for all $f\in\Mr$ and all $b\in \B$.

A $\rl$-algebra $\B$ is a \emph{$(\rl,\h)$-subalgebra} of
the \halg\ $\A$ if $\B$ is a $(\rl,\h)$-presubalgebra
of $\A$ such that $\io$ is an injective unital algebra
homomorphism.

In particular, this definition makes that a
part $\B_{\si,\tau}$ of the decomposition for $\B$
ends up in exactly one part $\A_{s(\si),s(\tau)}$
of the decomposition of $\A$. Requiring only the
weaker condition $\io(\B_{r(\al),r(\be)})\subseteq \A_{\al\be}$
does not imply this property. Note that in case
$\rl=\h$ we take $s$ to be the identity, and
we obtain the notions of $\h$-subprealgebra and
$\h$-subalgebra.

Using the map $s$ in the definition of
$\B$ being a $(\h,\rl)$-subprealgebra of $\A$
we extend the map $j\colon \Mr\to\Mh$ to
shift operators of the form
$fT_\si\in\Dr$ with $\si\in w(\B)$
by putting $j(gT_\si)= j(g)T_{s(\si)}\in\Dh$.
With this definition
$j\otimes \io\colon \Dr\wtt \B \to
\Dh\wtt \A$, which is indeed a well-defined
map, equals $\io\colon \B\to\A$ with the
identification $\Dr\wtt\B\cong\B$, $\Dh\wtt\A\cong\A$.
A similar remark applies to $\io\otimes j$.

\subsection{Singular vectors}\label{ssec:singularvectors}

We assume that $\U$ and $\A$ are paired \hbialg s and
that $\io\colon \W\to \U$ makes the $\rl$-prealgebra
$\W$ into a $(\rl,\h)$-subprealgebra of $\U$.

\begin{defn}\label{def:singularvectors}
{\rm (i)} Let $R\colon V\to V\wtt\A$
be a right corepresentation of $\A$ in the \hspa\ $V$.
A vector $v\in V$ is a $\W$-singular vector
if $\pi(\io(Y))v = 0$ for all $Y\in\W$
with $\pi(X) = (\id\otimes \pair{X}{\cdot}T_\be)\circ R$,
$X\in\U_{\al\be}$, as in Proposition
\ref{prop:cotodynrep}(i). \par\noindent
{\rm (ii)}
Let $L\colon V\to \A\wtt V$
be a left corepresentation of $\A$ in the \hspa\ $V$.
A vector $v\in V$ is a $\W$-singular vector
if $\pi(\io(Y))v = 0$ for all $Y\in\W$
with $\pi(X) = (T_\al\pair{X}{\cdot}\otimes\id)\circ L$,
$X\in\U_{\al\be}$, as in  Proposition
\ref{prop:cotodynrep}(ii).
\end{defn}

\begin{remark}\label{rmk:defsingularvectors}
(i) By decomposing a $\W$-singular
vector into homogeneous components $v=\sum_\ga v_\ga$,
$v_\ga\in V_\ga$, of the \hspa\ $V$, we obtain that
each of the homogeneous components is a
$\W$-singular vector. So we can assume without
loss of generality that the $\W$-singular
vector is homogeneous.
\par\noindent
(ii) Note that the space of $\W$-singular vectors
is a vector space over $\Mh$. For
the case of a right corepresentation $R$ we see
that for $Y\in
\W_{\si\tau}$ we have $\io(Y)\in\U_{s(\si),s(\tau)}$, so that for
$\W$-singular vector $v$ we have $\pi(\io(Y))(\mu_V(f)v)=
\mu_V(T_{-s(\si)}f)\, \pi(\io(Y))v = 0$.
\par\noindent
(iii) Assume that $R\colon V\to V\wtt \A$ is a right
corepresentation of $\A$ in the \hspa\ $V$, then
for $Y\in \W_{\si\tau}$ we have from
Proposition \ref{prop:cotodynrep}(i)
$\pi(\io(Y))V_\ga\subseteq V_{\ga+s(\si)-s(\tau)}$.
Hence $V_\ga$ consists of $\W$-singular vectors
in case the weight $\ga+s(\si)-s(\tau)$ does not
occur in $V$ for all $\si,\tau\in w(\W)$.
A similar remark applies to singular vectors
in left corepresentations.
\end{remark}

\begin{example}\label{eexample:singularvectorforE}
Let $\W$ be the \hprealg\ generated by $\be(z)$
viewed as subalgebra of $\E^\opp$. Consider the
left corepresentation $L$ of $\E$ on $V^N$
defined by \eqref{eq:defleftcorepU}.
It follows
from \eqref{eq:repAcopopp} that $v_N(z)$ is a $\W$-singular
vector.
Similarly, $v_0(z)$ is a $\W^\prime$-singular vector
for $\W^\prime$ be the \hprealg\ generated by $\ga(z)$
viewed as subalgebra of $\E^\opp$.
Combining this with the
terminology of Example \ref{ex:sphericalvectors},
$v_0(z)$
is lowest weight vector and $v_N(z)$
is highest weight vector with the weight
of $v_0(z)$, respectively $v_N(z)$, being $-N$
respectively $N$.
\end{example}

\subsection{Spherical vectors and spherical
corepresentations}\label{ssec:sphericalvectors}

We assume that $\U$ and $\A$ are paired \hbialg s, and
that $\io\colon \W\to \U$ makes the $\rl$-prealgebra
$\W$ into a $(\rl,\h)$-subprealgebra of $\U$.

\begin{defn}\label{def:sphericalvectors}
{\rm (i)} Let $R\colon V\to V\wtt\A$
be a right corepresentation of $\A$ in the \hspa\ $V$.
A vector $v\in V$ is a $\W$-spherical vector
if $\pi(\io(Y))v = \mu_V(\ep^\U(\io(Y))\1) v$ for all $Y\in\W$
with $\pi(X) = (\id\otimes \pair{X}{\cdot}T_\be)\circ R$,
$X\in\U_{\al\be}$, as in Proposition
\ref{prop:cotodynrep}(i). \par\noindent
{\rm (ii)}
Let $L\colon V\to \A\wtt V$
be a left corepresentation of $\A$ in the \hspa\ $V$.
A vector $v\in V$ is a $\W$-spherical vector
if $\pi(\io(Y))v = \mu_V(\ep^\U(\io(Y))\1) v$ for all $Y\in\W$
with $\pi(X) = (T_\al\pair{X}{\cdot}\otimes\id)\circ L$,
$X\in\U_{\al\be}$, as in Proposition
\ref{prop:cotodynrep}(ii).
\end{defn}

Note that Definition \ref{def:sphericalvectors}
coincides with Definition \ref{def:singularvectors}
in case $\ep^\U(\io(Y))\1=0$ for all $Y\in\W$.

\begin{remark}\label{rmk:defsphericalvectors}
(i) By decomposing a $\W$-spherical
vector into homogeneous components $v=\sum_\ga v_\ga$,
$v_\ga\in V_\ga$, of the \hspa\ $V$, we obtain that
each of the homogeneous components is a
$\W$-spherical vector. So we can assume without
loss of generality that the $\W$-spherical
vector is homogeneous. \par\noindent
(ii) Assume then that $v\in V_\ga$ is $\W$-spherical
for $V$ a right or left corepresentation of $\A$,
then for $g\in\Mr$ we have
$\mu_V(j(g))v = \mu_V(\ep^\U(\io(\mu_r^\W(g)))\1)v =
\pi(\io(\mu_r^\W(g)))v = \pi(\mu_r^\U(j(g)))v =
\mu_V(T_{-\ga}j(g))v$, see Proposition \ref{prop:cotodynrep}.
So $j(g)=T_{-\ga}j(g)$ for
all $g\in\Mr$, or $T_{-r(\ga)}g=g$, so we obtain $r(\ga)=0$.
\par\noindent (iii) Note that the space of $\W$-spherical vectors
is in general not a vector space over $\Mh$. For $Y\in
\W_{\si\tau}$ we have $\io(Y)\in\U_{s(\si),s(\tau)}$, so that for
$\W$-spherical vector $v$ in a right corepresentation
we have $\pi(\io(Y))(\mu_V(f)v)=
\mu_V(T_{-s(\si)}f)\, \pi(\io(Y))v = \mu_V(T_{-s(\si)}f)\, \ep^\U(\io(Y)\1)v$,
which in general is not equal to $\ep^\U(\io(Y)\1)(\mu_V(f)v)$ unless
$s(\si)=0$ or $\ep^\U(\io(Y)\1)=0$. For Example
\ref{ex:sphericalvectors} we are in the first case. Put
$\Mh(w(\W))$ to be the $\C$-linear
space of functions $f\in\Mh$
satisfying $T_{-s(\si)}f=f$ for all $\si\in w(\W)$ such that there
exists a $Y\in\W_{\si,\tau}$ for some $\tau\in\rls$ with
$\ep^\U(\io(Y)\1)\not=0$ in $\Mh$. Then we see that in general the
space of $\W$-spherical functions is a vector
space over $\Mh(w(\W))$, which is a space of meromorphic
functions on $\hs$ with specific periodicity properties.
\end{remark}

\begin{example}\label{ex:sphericalvectors}
Put $\Ih=\Mh\otimes_\C \Mh$,
and define the bigrading by $\Ih =(\Ih)_{00}$
and the left and right moment map
by $\mu_l^{\Ih}(f) =f\otimes \1$, $\mu^{\Ih}_r(g)=\1\otimes g$.
Then $\Ih$ is a \halg , see Rosengren \cite{Rose}
for more properties of $\Ih$.
For any \halg\ $\U$ the map $\io\colon \Ih\to
\U$, $\io(f\otimes g) = \mu_l^\U(f)\mu_r^\U(g)$
shows that $\Ih$ is a $\h$-subalgebra of $\U$.
With the notation of
Definition \ref{def:sphericalvectors} we find
that a spherical vector $v\in V$ in a
right corepresentation space  satisfies
$\mu_V(fg)v= \mu_V(\ep^\U(\io(f\otimes g)))v
= \pi(\mu_l^\U(f)\mu_r^\U(g))v =
\mu_V(f(T_{-\ga} g))v$ for all $f,g\in \Mh$,
so that $v$ is an $\Ih$-spherical vector
if and only if $v\in V_0$, and
$V_0$ is the space of $\Ih$-spherical
vectors. In this case $\Mh(w(\Ih))=\Mh$.

Note that we can extend this to
define the weight of a vector $v$ in the
representation space $V$ as $\ga\in\hs$ such that
$\pi(\io(f\otimes g)v =
\mu_V(f(T_{-\ga} g))v$, so this characterises $v\in V_\ga$.

In particular for the
left corepresentation $L$ of $\E$ on $V^N$
defined by \eqref{eq:defleftcorepU}, we see that
there is a spherical vector only for $N\in 2\N$, and
then $v_{N/2}(z)$ is the spherical vector.
\end{example}

We are interested in the matrix elements of a left or right
corepresentation corresponding to a spherical vector, and
in particular in their invariance properties.
Assume that
$\U$ and $\A$ are paired \hbialg s, that $R$ is a right
corepresentation of $\A$, and that $\W$ is a
$(\h,\rl)$-subprealgebra of $\U$.
Assume we have a basis $\{ v_i\}_{i\in I}$ of the \hspa\ $V$ such that $v_0$ is
a $\W$-spherical vector, and let $Rv_i=\sum_{k\in I} v_k\otimes
R_{ki}$. Using Lemma \ref{lemma:interactioncorep-rep}(i) we find
\begin{equation*}
\begin{split}
\sum_{k\in I} v_k \otimes \io(Y)\cdot R_{k0}
&= R\bigl( \pi(\io(Y))v_0 \bigr)
= R\bigl( \mu_V(\ep^\U(\io(Y))\1)v_0\bigr) \\
&= \mu_{V\wtt\A} (\ep^\U(\io(Y))\1) Rv_0
= \sum_{k\in I} v_k \otimes  \mu_r^\A(\ep^\U(\io(Y))\1)R_{k0},
\end{split}
\end{equation*}
which implies
$\io(Y)\cdot R_{k0}=\mu_r^\A(\ep^\U(\io(Y))\1)R_{k0}$
for all $Y\in\W$.

\begin{prop}\label{prop:invpropsphericalfunction}
Assume that $\U$ and
$\A$ are paired \hHopfster s,
and that $\W$ is a $(\h,\rl)$-subprealgebra of $\U$
given by $\io\colon\W\to\U$. \par\noindent
{\rm (i)} Let $\{ v_i\}_{i\in I}$ be a homogeneous
orthogonal basis for the right unitarisable
corepresentation $R\colon V\to
V\wtt \A$, $Rv_i=\sum_{k\in I} v_k\otimes R_{ki}$
such that $\pair{v_i}{v_j}=\de_{ij}N_j$.
Assume that $v_0\in V$ is a $\W$-spherical vector, then
\begin{equation}\label{eq:sphericalfunctionsinR}
\begin{split}
&\io(Y)\cdot (\mu_l^\A(N_0)R_{00}) =\mu_r^\A(\ep^\U(\io(Y))\1)(\mu_l^\A(N_0)R_{00}), \\
& (\mu_l^\A(N_0) R_{00})\cdot \io(Y)^\ast = \mu_l^\A(\ep^\U(\io(Y)^\ast)\1)
(\mu_l^\A(N_0) R_{00}).
\end{split}
\end{equation}
{\rm (ii)} Let $\{ v_i\}_{i\in I}$ be a homogeneous
orthogonal basis for the left unitarisable
corepresentation $L\colon V\to
\A\wtt V$, $Lv_i=\sum_{k\in I} L_{ik} \otimes v_k$
such that $\pair{v_i}{v_j}=\de_{ij}N_j$.
Assume that $v_0\in V$ is a $\W$-spherical vector, then
\begin{equation}\label{eq:sphericalfunctionsinL}
\begin{split}
& (\mu_r^\A(N_0) L_{00})\cdot \io(Y) = \mu_l^\A(\ep^\U(\io(Y))\1)
(\mu_r^\A(N_0) L_{00}),\\
&\io(Y)^\da \cdot (\mu_r^\A(N_0)L_{00}) =\mu_r^\A(\ep^\U(\io(Y)^\da)\1)
(\mu_r^\A(N_0)L_{00}).
\end{split}
\end{equation}
\end{prop}

\begin{proof} For the first statement, the first equality
has been derived in more generality and
\eqref{eq:actionmomentA} gives the first equation of
\eqref{eq:sphericalfunctionsinR}.
For the second we use the first
in combination with \eqref{eq:propactiononmatrixeltsunitarycorepR}
to find
\begin{equation}
\begin{split}
\bigl(\mu_l^\A(N_0) R_{0j}\bigr) \cdot \io(Y)^\ast &=
S^\A\bigl( \io(Y)\cdot (\mu_l^\A(N_j) R_{j0})\bigr)^\ast =
S^\A\bigl( \mu_r^\A(\ep^\U(\io(Y)\1)) \mu_l^\A(N_j) R_{j0})\bigr)^\ast \\ &=
\mu_r^\A(N_j) \mu_l^\A(\ep^\U(\io(Y)^\ast\1)) \, S^\A(R_{j0})^\ast \\
&=
\mu_r^\A(N_j) \mu_l^\A(\ep^\U(\io(Y)^\ast\1)) \,
\mu_l^\A(N_0)\mu_r^\A(N_j^{-1}) R_{0j} \\ &=
\mu_l^\A(\ep^\U(\io(Y)^\ast\1)) \,
\mu_l^\A(N_0) R_{0j}
\end{split}
\end{equation}
using \eqref{eq:propactiononmatrixeltsunitarycorepR} again
and taking $j=0$ gives the result. The
second statement follows similarly.
\end{proof}

\begin{example}\label{ex2:sphericalvectors}
In the case $\W=\Ih$ as in
Example \ref{ex:sphericalvectors},
Proposition \ref{prop:invpropsphericalfunction}
says that the spherical functions
$\mu_l^\A(N_0)R_{00}$ and
$\mu_r^\A(N_0) L_{00}$  are elements of $\A_{00}$.
\end{example}

Because of Proposition \ref{prop:invpropsphericalfunction}
it is natural to consider the spaces
\begin{equation}\label{eq:spaceofsphericalfR}
\begin{split}
\A^\W = \{ a \in \A \mid
\io(Y)\cdot a =\mu_r^\A(\ep^\U(\io(Y))\1)a, \quad
 a\cdot \io(Y)^\ast = \mu_l^\A(\ep^\U(\io(Y)^\ast)\1)a,\
 \forall Y\in\W \}, \\
{}^\W\A = \{ a \in \A \mid
\io(Y)^\da \cdot a =\mu_r^\A(\ep^\U(\io(Y)^\da)\1)a, \quad
 a\cdot \io(Y) = \mu_l^\A(\ep^\U(\io(Y))\1)a,\
 \forall Y\in\W \}
\end{split}
\end{equation}
for a right or left unitarisable corepresentation,
but in general we need more information on $\W$ in order
to be able to say more on $\A^\W$ and ${}^\W\A$,
except being $\C$-linear spaces.
E.g. in case $\ep^\U\circ\io \colon \W\to \Dh$
is zero, we see that $\A^\W$ and ${}^\W\A$ are
invariant under multiplication by $\mu_l^\A(f)$,
$\mu_r^\A(f)$ for any $f\in\Mh$, and if moreover
$\De^\U(\io(\W))\subset \U\wtt \io(\W) + \io(\W)\wtt\U$
we see that $\A^\W$ and ${}^\W\A$ are $\h$-subalgebras
of $\A$ by Theorem \ref{thm:defaction}.

So in this context it is natural to
use  the notion of $\h$-coideal
as introduced by Rosengren \cite[Def.~4.5]{Rose}.

\begin{defn}\label{def:coideal} A $\C$-linear subspace $I$ of a
\hbialg\ $\U$ is a $\h$-coideal if
\begin{enumerate}
\item[(i)] $I$ is a 2-sided ideal in $\U$ as an
associative algebra,
\item[(ii)] $(\U_{\al\be}+I)\cap (\U_{\ga\de}+I) = I$
for $(\al,\be)\not= (\ga,\de)$,
\item[(iii)] $\De^\U(X) \in \U\wtt I + I\wtt\U$
for all $X\in I$,
\item[(iv)] $\ep^\U(X)=0$ for all $X\in I$.
\end{enumerate}
\end{defn}

Assuming that $\U$ and $\A$ are paired \hbialg s
and that $I\subset \U$ is a $\h$-coideal, we have the
following Proposition, whose proof is left to the reader.

\begin{prop}\label{prop:sphericalsubspace}
Let $\U$ and $\A$ be paired \hbialg s
and assume that $I\subset \U$ is a $\h$-coideal.
Put $\A^I = \{ a \in A \mid X\cdot a = 0 = a\cdot X,
\ \forall \, X\in I\}$,
then $\A^I$ is a $\h$-subalgebra of $\A$.
Moreover, in case $\U$ and $\A$ are paired
\hHopfster s and $I$ is $\ast\circ S^\U$-invariant,
then $\A^I$ is  $\ast$-invariant.
\end{prop}

\begin{example}\label{ex:hcoideal} In order to give an
example of a $\h$-coideal, we first consider group-like elements in
a \hcoalg\ or in a \hbialg . We say that a
homogeneous element $X\in\U_{\al\be}$ is
a \emph{group-like element} if $\De^\U(X)=X\otimes X$
and $\ep^\U(X)=T_{-\al}$. Note that this forces
$\be=\al$. For two group-like elements $X,Y\in\U_{\al\al}$
we have
\begin{equation*}
\De^\U(X-Y) = (X-Y)\otimes X + Y\otimes (X-Y),\qquad
\ep^\U(X-Y)=0,
\end{equation*}
so the $2$-sided ideal generated by differences of group like
elements satisfies Definition \ref{def:coideal}(i), (iii), (iv),
and Definition \ref{def:coideal}(ii) has to be checked
in specific examples.

Consider an invertible element $f\in\Mh$, and let
$\mu_l^\U(f) \mu_r^\U(f^{-1}) \in \U_{00}$
then
$\ep^\U\bigl(\mu_l^\U(f)\mu_r^\U(f^{-1}) \bigr)=
f f^{-1}=1\in\Dh$ and
\begin{equation*}
\De^\U\bigl( \mu_l^\U(f)\mu_r^\U(f^{-1}) \bigr)
= (\mu_l^\U(f))\otimes (\mu_r^\U(f^{-1}))
= (\mu_l^\U(f)\mu_r^\U(f^{-1}))
\otimes (\mu_l^\U(f)\mu_r^\U(f^{-1}))
\end{equation*}
using \eqref{eq:defotimesMh} in the last equality.
So $\mu_l^\U(f) \mu_r^\U(f^{-1})$ is a group-like
element for any invertible $f\in\Mh$. Let
$I$ be the $2$-sided ideal generated by
$\mu_l^\U(f) \mu_r^\U(f^{-1}) - \mu_l^\U(g) \mu_r^\U(g^{-1})$,
$f,g\in\Mh$ invertible. Then it is easily checked that
$I$ is a $\h$-Hopf coideal, which is $S^\U$-invariant in case
$\U$ is a \hHopfalg\ and moreover $\ast$-invariant if
$\U$ is \hHopfster .
In case $\U$ and $\A$ are paired \hHopfalg s, we have
$\A^I=\A_{00}$.
\end{example}



\end{document}